\documentclass{article}%
\usepackage{amsfonts}
\usepackage{amsmath}
\usepackage{amssymb}
\usepackage{amsxtra}
\usepackage{graphicx}
\usepackage{geometry}
\usepackage{color}
\usepackage{colortbl}%
\providecommand{\U}[1]{\protect\rule{.1in}{.1in}}
\newtheorem{theorem}{Theorem}

\newtheorem{lemma}[theorem]{Lemma}

\newtheorem{remark}[theorem]{Remark}

\begin{document}

\title{Asymptotic analysis and numerical modeling of mass transport in tubular structures}
\author{G.Cardone\\Universit\`{a} del Sannio - Dipartimento di Ingegneria\\Piazza Roma, 21 - 84100 Benevento, Italy\\email: giuseppe.cardone@unisannio.it
\and G.P.Panasenko\\University Jean Monnet - LaMUSE \\Saint Etienne, 23, rue P.Michelon, 42023 St. Etienne, France\\email: grigory.panasenko@univ-st-etienne.fr
\and Y.Sirakov\\University Jean Monnet - LaMUSE \\Saint Etienne, 23, rue P.Michelon, 42023 St. Etienne, France\\email: Yvan.Sirakov@univ-st-etienne.fr}
\maketitle

\begin{abstract}
In the paper the flow in a thin tubular structure is considered. The velocity
of the flow stands for a coefficient in the diffusion-convection equation set
in the thin structure. An asymptotic expansion of solution is constructed.
This expansion is used further for justification of an asymptotic domain
decomposition strategy essentially reducing the memory and the time of the
code. A numerical solution obtained by this strategy is compared to the
numerical solution obtained by a direct FEM computation.

\medskip



\end{abstract}

\section{Introduction}

The paper is devoted to the strategy of numerical implementation of the
asymptotic partial decomposition of the domain for the tubular structures of a
complicated geometry. We will consider the Stokes flow in this structure and
the convection-diffusion and sorption process for some diluted substance.
First we consider an asymptotic expansion of the solution. We emphasize the
importance of the boundary layers in the neighborhood of some special
structural elements of the tubular domain, such as the bifurcations of canals
and "stenosis areas". That is why some multiscale strategy should be applied
to the analysis of the convection-diffusion process: the 1D limit description
in the canals will be coupled with some 2D zooms in these special structural elements.

In section 2 we define a tubular structure as a union of thin rectangles
connected by some domains of small diameters. The Stokes equation and the
diffusion-convection equation are set in this domain. For the Stokes equation
the Dirichlet conditions are respected at the lateral boundary with some given
inflow and outflow. For the diffusion-convection equation we pose the Robin
type condition at the lateral boundary with some given inflow and outflow
concentrations. The viscosity and the diffusion are constant out of some
"stenosis area" where they may have variations. The varying viscosity can be
used for the modeling of a clot in the blood circulation process. Indeed, if
at some part of the domain the viscosity is great, then, applying the idea of
the fictitious domain method, we can exclude this part of the domain from the
flow area (see Remark \ref{rem1}).

In section 3 we consider the Stokes equation in tubular structure. The
asymptotic expansion of the solution for constant viscosity has been obtained
in \cite{GPbook}. In this section we construct the boundary layer correctors
for the varying viscosity in the stenosis areas.

In section 4 the convection-diffusion equation is considered. First we
construct the asymptotic expansion in an infinite tube (subsection 4.1). Then
using this expansion as a regular ansatz we add the boundary layer correctors
in the stenosis zones (subsection 4.2) in a bifurcation area and in the
entrance/exit elements (subsections 4.3 and 4.4). In the subsection 4.5 the
leading term of the asymptotic expansion is presented. The justification of
the asymptotic expansion follows the scheme: estimate for the residuals and
application of the a priori estimates for the initial problem.

Section 5 describes one version of the partial asymptotic domain decomposition
strategy for the mass transport problem in a tubular structure.

Finally, section 6 develops the numerical experiment comparing the direct
numerical solution of the 2D problem and the asymptotic solution of the
partially decomposed problem. The results of this experiment confirm good
coincidence of the exact solution and the approximate solution obtained by the
method of asymptotic partial decomposition of domain.

\section{Geometry of tubular structure and setting of the problem}

We will introduce the tubular domain which consists of three types of
structural elements: canals, bifurcations and stenosis areas. This tubular
structure is similar to the rod structures introduced in \cite{GP91} and the
tube structures or pipe structures introduced in \cite{GP00}; we consider a
new element that is, the stenosis area, simulated by varying coefficients of
the equation (viscosity and diffusion coefficients) and not by geometric singularity.

Let us remind the definition of a tube structure.

Let $e_{1},...,e_{n}$ be n closed segments in $\mathbb{R}^{2}$ which have a
single common point $0$ (i.e. the origin of the coordinate system) and let it
be the common end point of all these segments. Let $\theta_{1},...,\theta
_{n}\in\left(  0,1\right)  $ be $n$ positive numbers. Making a change of
variables (rotation) such that the new axis $x_{1}$ denoted $x_{1}^{e_{i}}$
contains the segment $e_{i}$ and the second new axis $x_{2}^{e_{i}}$ is
orthogonal to $e_{i}$, we define
\[
B_{i}^{\varepsilon}=\left\{  \left(  x_{1},x_{2}\right)  \in\mathbb{R}%
^{2}:x_{1}^{e_{i}}\in\left(  0,|e_{i}|\right)  ,\ x_{2}^{e_{i}}\in\left(
-\frac{\theta_{i}\varepsilon}{2},\frac{\theta_{i}\varepsilon}{2}\right)
\right\}  .
\]%
\[%
%TCIMACRO{\FRAME{itbpFU}{2.2477in}{1.3733in}{0in}{\Qcb{Fig. 1}}{}%
%{fig1.eps}{\special{ language "Scientific Word";  type "GRAPHIC";
%width 2.2477in;  height 1.3733in;  depth 0in;  original-width 2.2087in;
%original-height 1.3387in;  cropleft "0";  croptop "1";  cropright "1";
%cropbottom "0";  filename 'fig1.eps';file-properties "XNPEU";}} }%
%BeginExpansion
{\parbox[b]{2.2477in}{\begin{center}
\includegraphics[
height=1.3733in,
width=2.2477in
]%
{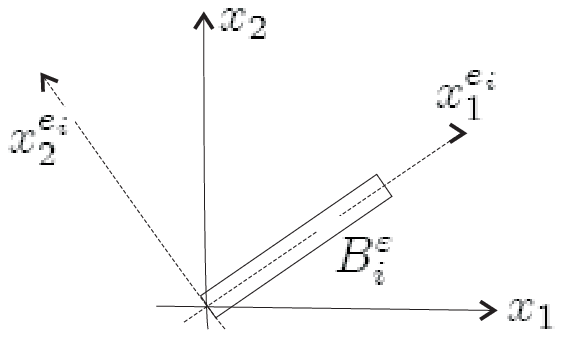}%
\\
Fig. 1
\end{center}}}
%EndExpansion
\]
Let $\gamma_{0}$ be a bounded domain containing $O$, $\gamma_{i}$ be a bounded
domain containing $O_{i},$ the end point of $e_{i}$ (different from $O$). We
assume for simplicity that diam$\left(  \gamma_{i}\right)  ,$ diam$\left(
\gamma_{O}\right)  <2.$ Consider the homothetic contraction of $\gamma_{0}$ in
$\dfrac{1}{\varepsilon}$ times with the center of the homothety in $O$ and
denote $\gamma_{O}^{\varepsilon}$ the image of $\gamma_{0}$ by this homothety,
i.e.%
\[
\gamma_{O}^{\varepsilon}=\left\{  \left(  x_{1},x_{2}\right)  \in
\mathbb{R}^{2}:\left(  \frac{x_{1}}{\varepsilon},\frac{x_{2}}{\varepsilon
}\right)  \in\gamma_{0}\right\}  .
\]

In the same way we consider%
\[
\gamma_{i}^{\varepsilon}=\left\{  \left(  x_{1},x_{2}\right)  \in
\mathbb{R}^{2}:\frac{\left(  x_{1},x_{2}\right)  -O_{i}}{\varepsilon}+O_{i}%
\in\gamma_{i}\right\}
\]
the homothetic contraction of $\gamma_{i}$ in $\dfrac{1}{\varepsilon}$ times
with the center of the homothety in $O_{i}$.

Define the one bundle tube structure%
\[
B_{\varepsilon}=\left(
%TCIMACRO{\tbigcup \limits_{i=1}^{n}}%
%BeginExpansion
{\textstyle\bigcup\limits_{i=1}^{n}}
%EndExpansion
B_{i}^{\varepsilon}\right)  \cup\left(
%TCIMACRO{\tbigcup \limits_{i=1}^{n}}%
%BeginExpansion
{\textstyle\bigcup\limits_{i=1}^{n}}
%EndExpansion
\gamma_{i}^{\varepsilon}\right)
\]
such that it is a connected domain with the $C^{2}-$smooth boundary.

In a more general case, we will consider several one-bundle structures
$B_{\varepsilon1},B_{\varepsilon2},...,B_{\varepsilon m}$ such that every of
these structures is associated to some segments:%
\begin{align*}
&  e_{11},...,e_{1n_{1}}\text{ \ for }B_{\varepsilon1},\\
&  e_{21},...,e_{2n_{2}}\text{ \ for }B_{\varepsilon2},\\
&  ...\\
&  e_{m1},...,e_{mn_{n}}\text{ \ for }B_{\varepsilon m}.
\end{align*}
Assume that if two of these segments have common point then it is an end point
for each of these segments; assume that the union $\bigcup\limits_{q=1}%
^{m}B_{\varepsilon q}$ is a connected domain with $C^{2}-$smooth boundary. In
this case we will call this union $B_{\varepsilon}=\bigcup\limits_{q=1}%
^{m}B_{\varepsilon q}$ a multi-bundle tube structure. If an end point of some
segment $e_{ij}$ is not an end point for all other segments then such end
point will be called solitary.

We will consider the Stokes equation and the convection-diffusion equation in
such tube structure. The boundary condition is the vanishing velocity for the
Stokes equation everywhere except of some special parts of the boundary
(entrance and exit). These parts are some connected parts $\Gamma_{i}$ of the
boundary of smoothing domains $\gamma_{i}.$ We assume that the end points of
segments corresponding to this $\gamma_{i}$ are solitary end points. Let
$\Gamma_{1},...,\Gamma_{r}$ be these parts of the boundary.

\bigskip

We consider the Stokes equation in such a tube structure with the varying
viscosity coefficient $\mu:$%

\begin{equation}%
\begin{array}
[c]{l}%
\operatorname{div}\left(  \mu_{\varepsilon}\left(  x\right)  \left(  \nabla
u_{\varepsilon}+\left(  \nabla u_{\varepsilon}\right)  ^{T}\right)
-p_{\varepsilon}I\right)  =f_{\varepsilon}\left(  x\right)  ,\\
\operatorname{div}u_{\varepsilon}=0,
\end{array}
\label{problin}%
\end{equation}
where the divergence is taken with respect to the elements of each line of the
matrix $\mu_{\varepsilon}\left(  x\right)  \left(  \nabla u_{\varepsilon
}+\left(  \nabla u_{\varepsilon}\right)  ^{T}\right)  -p_{\varepsilon}I,$ and
the convection-diffusion equation%
\begin{equation}
-\operatorname{div}\left(  K_{\varepsilon}\left(  x\right)  \nabla
c_{\varepsilon}\right)  +u_{\varepsilon}\left(  x\right)  \cdot\nabla
c_{\varepsilon}=g\left(  x_{1}^{e_{i}}\right)  . \label{1.1000}%
\end{equation}
Assume that $g=0$ in some neighborhood of the end points of the segments,
$g\in C^{k+2}\left(  e_{i}\right)  $ for all segments $e_{i}.$

We will assume that $\mu_{\varepsilon}$ and $K_{\varepsilon}$ are positive
constants $\mu$ and $\varkappa$ respectively, everywhere except some "stenosis
areas" where they have a form:%
\begin{equation}
\mu_{\varepsilon}\left(  x\right)  =\mu+M\left(  \frac{x-\overline{x}_{s}%
}{\varepsilon}\right)  \text{ \ and \ }K_{\varepsilon}\left(  x\right)
=\varkappa+K\left(  \frac{x-\overline{x}_{s}}{\varepsilon}\right)
\label{1.600}%
\end{equation}
where $M$ and $K$ are measurable bounded function having a finite support
inside the ball $B\left(  \dfrac{\overline{x}_{s}}{\varepsilon},2\right)  $
with the center $\dfrac{\overline{x}_{s}}{\varepsilon}$ and the radius $2,$
such that,
\[
\exists\varkappa_{1}>0:\mu_{\varepsilon}\left(  x\right)  ,\ K_{\varepsilon
}\left(  x\right)  \geq\varkappa_{1}.
\]
Here $\overline{x}_{s}$ are some points belonging to the segments $e_{i}$ of
the graph of the structure, they are different from the end points and are
independent of $\varepsilon.$

The sorption will be modeled by the boundary condition for the
diffusion-convection equation, i.e. let us consider the boundary conditions:%
\begin{align}
u_{\varepsilon}  &  =0\text{ on the lateral boundary }\partial B_{\varepsilon
}\setminus\left(
%TCIMACRO{\tbigcup \limits_{t=1}^{r}}%
%BeginExpansion
{\textstyle\bigcup\limits_{t=1}^{r}}
%EndExpansion
\Gamma_{t}\right)  ,\label{1.2010}\\
K_{\varepsilon}\left(  x\right)  \frac{\partial c_{\varepsilon}}{\partial n}
&  =\varepsilon\beta c_{\varepsilon}\text{ \ on the lateral boundary }\partial
B_{\varepsilon}\setminus\left(
%TCIMACRO{\tbigcup \limits_{t=1}^{r}}%
%BeginExpansion
{\textstyle\bigcup\limits_{t=1}^{r}}
%EndExpansion
\Gamma_{t}\right)  ,\label{1.2020}\\
u_{\varepsilon}  &  =G\left(  \frac{x-x_{b_{t}}}{\varepsilon}\right)  \text{
\ on }\Gamma_{t}\label{1.2030}\\
c_{\varepsilon}  &  =q_{t}=const\text{ \ on }\Gamma_{t},\text{ }t=1,...,r
\label{1.2040}%
\end{align}
where $x_{b_{t}}$ is an end point, inside $\gamma_{t},$ of a corresponding
segment, $n$ is an outer normal, $G\in C_{0}^{2}\left(  \overline{\gamma}%
_{t}\right)  $ and $\sum_{t}%
%TCIMACRO{\dint _{\Gamma_{t}}}%
%BeginExpansion
{\displaystyle\int_{\Gamma_{t}}}
%EndExpansion
n\cdot G\left(  \frac{x-x_{b_{t}}}{\varepsilon}\right)  ds=0.$

\begin{remark}
\label{rem1}The varying viscosity and diffusion coefficients can be used for
the modelling of a clot in the blood circulation process. Let us remind the
fictitious domain method. Consider an example: the Poisson equation $\Delta
u=f$ posed in a bounded domain $G$ with the boundary condition $u|_{\partial
G}=0;$ $f\in L^{2}\left(  G\right)  .$ The fictitious domain method reduces
this problem to the problem set in a larger rectangular $R\supset G:$%
\begin{align*}
\operatorname{div}\left(  K_{\omega}\left(  x\right)  \nabla u_{\omega
}\right)   &  =F\left(  x\right)  ,\ \ x\in R,\\
u_{\omega}|_{\partial R}  &  =0
\end{align*}
where%
\[%
\begin{array}
[c]{cc}%
K_{\omega}\left(  x\right)  =\left\{
\begin{array}
[c]{l}%
1,\ x\in G,\\
\omega,\ x\in R\diagdown G,
\end{array}
\right.  & F\left(  x\right)  =\left\{
\begin{array}
[c]{l}%
f\left(  x\right)  ,\ x\in G,\\
0,\ x\in R\diagdown G.
\end{array}
\right.
\end{array}
\]%
\[%
%TCIMACRO{\FRAME{itbpFU}{0.8415in}{0.8389in}{0in}{\Qcb{Fig. 2}}{}%
%{figure2.eps}{\special{ language "Scientific Word";  type "GRAPHIC";
%width 0.8415in;  height 0.8389in;  depth 0in;  original-width 1.3898in;
%original-height 1.3846in;  cropleft "0";  croptop "1";  cropright "1";
%cropbottom "0";  filename 'Figure2.eps';file-properties "XNPEU";}} }%
%BeginExpansion
{\parbox[b]{0.8415in}{\begin{center}
\includegraphics[
height=0.8389in,
width=0.8415in
]%
{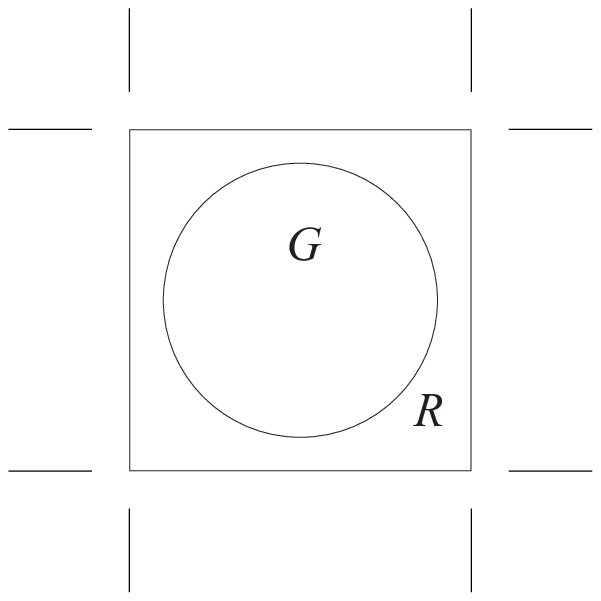}%
\\
Fig. 2
\end{center}}}
%EndExpansion
\]
For smooth $\partial G,$ one can prove that, as $\omega\rightarrow+\infty$
\[
u_{\omega}\rightarrow\left\{
\begin{array}
[c]{l}%
u\left(  x\right)  ,\ x\in G,\\
0,\ x\in R\diagdown G,
\end{array}
\right.  \ \ \text{in }H^{1}\left(  R\right)  .
\]
So the field $u_{\omega}$ vanishes in the fictitious part of the domain. The
same effect holds for the Stokes equation, where the varying viscosity can be
used to modelling the absence of flow in the fictitious part of domain
occupied by the clot.
\end{remark}

\section{The Stokes equation}

We will consider separately the problem for the Stokes equation and the
convection-diffusion one. First we apply the results of \cite{GP00} and get
the asymptotic expansion of the solution. At the second stage we assume that
the velocity $u_{\varepsilon}$ is known and consider the convection-diffusion
equation with the velocity coefficient corresponding to the first term of the
asymptotic approximation.

The asymptotic solution of the Stokes problem was considered in \cite{GPbook}.
The only difference is related to the "stenosis areas" where the boundary
layers are constructed as follows.

A stenosis area can be simulated by a varying viscosity in some close
neighborhood of the origin of the coordinate system. Then we can consider one
channel parallel to the $Ox_{1}$ axis, i.e. $\widehat{G}_{\varepsilon}=\left(
-1,1\right)  \times\left(  -\dfrac{\varepsilon}{2},\dfrac{\varepsilon}%
{2}\right)  .$ We will not take care of the ends of this channel because these
ends are supposed to be junction points with some other channels and the
construction of the bifurcation boundary layers is described in \cite{GPbook}.
So we will try to construct a solution of problem (\ref{problin}) stabilizing
to a Poiseuille solution as $\dfrac{x_{1}}{\varepsilon}\rightarrow+\infty.$

Let equations (\ref{problin}) be considered in this channel $\widehat
{G}_{\varepsilon}$ and let the viscosity coefficient $\mu_{\varepsilon}$ have
a structure%
\[
\mu_{\varepsilon}\left(  x\right)  =\mu+M\left(  \frac{x}{\varepsilon}\right)
\]
where $\mu>0$ is constant and $M$ has a support inside the ball $B\left(
0,2\right)  =\left\{  \xi\in\mathbb{R}^{2}:\xi_{1}^{2}+\xi_{2}^{2}<4\right\}
$ such that $M\left(  \xi\right)  \geq0.$

Then if the right-hand side $f_{\varepsilon}$ is equal to zero, the asymptotic
solution out of the boundary layer zone (at some finite distance from $0)$ is
a Poiseuille flow:%
\begin{equation}
\overline{u}_{p}=c_{1}\binom{\frac{1}{2\mu}\left(  x_{2}^{2}-\frac
{\varepsilon^{2}}{4}\right)  }{0},\ \ \ \overline{p}\left(  x\right)
=c_{1}x_{1}+c_{2}, \label{1.3000}%
\end{equation}
where $c_{1},c_{2}$ are some constants.

Then the boundary layer corrector has a form $\left(  \varepsilon^{2}U\left(
\xi\right)  ,\varepsilon P\left(  \xi\right)  \right)  $ and $\left(
U,P\right)  $ is a solution of the following problem%
\[%
\begin{array}
[c]{l}%
\operatorname{div}_{\xi}\left(  \left(  \mu+M\left(  \frac{x}{\varepsilon
}\right)  \right)  \left(  \nabla_{\xi}U+\left(  \nabla_{\xi}U\right)
^{T}\right)  -P\left(  \xi\right)  I\right)  =-c_{1}\operatorname{div}_{\xi
}\left(  M\left(  \xi\right)  \left(
\begin{array}
[c]{cc}%
0 & \frac{1}{\mu}\\
\frac{1}{\mu} & 0
\end{array}
\right)  \xi_{2}\right)  ,\\
\operatorname{div}_{\xi}U=0,\ \ \xi\in\left(  -\infty,+\infty\right)
\times\left(  -\dfrac{1}{2},\dfrac{1}{2}\right)  ,\\
U|_{\xi_{2}=\pm\frac{1}{2}}=0,
\end{array}
\]
where the right hand side
\[
-c_{1}\operatorname{div}_{\xi}\left(  M\left(  \xi\right)  \left(
\begin{array}
[c]{cc}%
0 & \frac{1}{\mu}\\
\frac{1}{\mu} & 0
\end{array}
\right)  \xi_{2}\right)
\]
has a support inside $B\left(  0,2\right)  .$ It is well known (\cite{Naz,
Galdi}) that this problem has a unique solution $\left(  U,P\right)  $
stabilizing to $\left(  0,const\right)  $ at the infinity. Then $\left(
\overline{u}_{p}+\varepsilon^{2}U\left(  \frac{x}{\varepsilon}\right)
,\overline{p}\left(  x\right)  +\varepsilon P\left(  \frac{x}{\varepsilon
}\right)  \right)  $ satisfies equation (\ref{problin}) and flow coincides
with the Poiseuille flow at a finite distance from zero with an exponentially
small error $O\left(  e^{-\frac{\alpha}{\varepsilon}}\right)  ,$ $\alpha>0.$
If $P\left(  \xi\right)  \rightarrow0$ for $\xi_{1}\rightarrow-\infty$ and
$P\left(  \xi\right)  \rightarrow c_{+}$ for $\xi_{1}\rightarrow+\infty,$ then
we have to proceed the following special gluing of pressure at zero:

1) we redefine a new (discontinuous)%
\[
\overline{p}\left(  x\right)  =\left\{
\begin{array}
[c]{l}%
c_{1}x_{1}+c_{2},\ \ \text{for }x_{1}<0,\\
c_{1}x_{1}+c_{2}+\varepsilon c_{+},\ \ \text{for }x_{1}>0;
\end{array}
\right.
\]

2) we redefine a new (discontinuous) (denoted by the same letter):%
\[
P\left(  \xi\right)  :=\left\{
\begin{array}
[c]{l}%
P\left(  \xi\right)  ,\ \ \ \ \ \ \ \ \ \ \text{for }\xi_{1}<0,\\
P\left(  \xi\right)  -c_{+},\ \ \text{for }\xi_{1}>0.
\end{array}
\right.
\]
Mention that the sum $\overline{p}\left(  x\right)  +\varepsilon P\left(
\frac{x}{\varepsilon}\right)  $ is still smooth and this new $P\left(
\xi\right)  \rightarrow0$ for $\xi_{1}\rightarrow\pm\infty,$ i.e. it has a
standard boundary layer shape. The justification of this expansion follows
laterally \cite{GPbook}

\section{The convection-diffusion equation}

Consider the diffusion-convection problem (\ref{1.1000}), (\ref{1.2020}) and
(\ref{1.2040}) in a tube structure $B_{\varepsilon},$ where the coefficient
$K_{\varepsilon}$ is given by formula (\ref{1.600}) and $u_{\varepsilon}$ in
(\ref{1.1000}) is replaced by the given vector-valued function $V_{\varepsilon
}$ having the following structure:%
\[
V_{\varepsilon}\left(  x\right)  =\left\{
\begin{array}
[c]{l}%
V_{i}\left(  \frac{x-O_{i}}{\varepsilon}\right)  \text{ \ \ for }x\in
B_{\varepsilon}:\left\vert \frac{x-O_{i}}{\varepsilon}\right\vert
<2,\ i=0,..,n,\\
\widetilde{V}_{s}\left(  \frac{x-\overline{x}_{s}}{\varepsilon}\right)  \text{
\ \ for }x\in B_{\varepsilon}:\left\vert \frac{x-\overline{x}_{s}}%
{\varepsilon}\right\vert <2,\\
\binom{V_{p}\left(  \frac{x_{2}^{e_{i}}}{\varepsilon}\right)  }{0}%
\ \ \text{for }x\in B_{\varepsilon}:x_{2}^{e_{i}}\in e_{j},\ \left\vert
\frac{x-O_{i}}{\varepsilon}\right\vert \geq2,\ \left\vert \frac{x-\overline
{x}_{s}}{\varepsilon}\right\vert \geq2,\text{ for all }i\text{ and }s,
\end{array}
\right.
\]
where $V_{p}\left(  \xi_{2}^{e_{i}}\right)  =-\left(  \left(  \xi_{2}^{e_{i}%
}\right)  ^{2}-\frac{1}{4}\right)  $ and $\widetilde{V}_{s}\left(  \xi\right)
$ and $V_{i}\left(  \xi\right)  $ are some given smooth vector-valued
functions with finite support in the ball $B\left(  0,2\right)  ,$ and
$\overline{x}_{s}$ are the "stenosis nodes" and $O_{i}$ are the ends of the
segments $e_{i}.$

Really, the structure of the velocity field $u_{\varepsilon}$ is more
complicated: out of the balls of radius $2\varepsilon$ surrounding the nodes
$O_{i}$ and the stenosis points $\overline{x}_{s}$ the velocity differs from
the Poiseuille flow by some exponentially decaying boundary layer functions.
Here we simplify the structure of the velocity field replacing $u_{\varepsilon
}$ by $V_{\varepsilon}.$ So we consider here the problem:%
\begin{equation}
-\operatorname{div}\left(  K_{\varepsilon}\left(  x\right)  \nabla
c_{\varepsilon}\right)  +V_{\varepsilon}\left(  x\right)  \cdot\nabla
c_{\varepsilon}=g\left(  x_{1}^{e_{i}}\right)  \label{2.000}%
\end{equation}
with the boundary conditions%
\begin{align}
K_{\varepsilon}\left(  x\right)  \frac{\partial c_{\varepsilon}}{\partial n}
&  =\varepsilon\beta c_{\varepsilon}\text{ \ on the lateral boundary }\partial
B_{\varepsilon}\setminus\left(
%TCIMACRO{\tbigcup \limits_{t=1}^{r}}%
%BeginExpansion
{\textstyle\bigcup\limits_{t=1}^{r}}
%EndExpansion
\Gamma_{t}\right)  ,\label{2.100}\\
c_{\varepsilon}  &  =q_{t}=const\text{ \ on }\Gamma_{t},\text{ }t=1,...,r.
\label{2.200}%
\end{align}

We assume that $\beta$ is a constant and $\beta\leq0.$ It means that the
sorption takes place: the outflow $-\varkappa\dfrac{\partial c_{\varepsilon
r}}{\partial n}$ is positive. It is well known that if at least one of given
concentrations $c_{s}$ is equal to zero, then the Poincar\'{e}-Friedrichs
inequality for $c_{\varepsilon}$ holds with a constant independent of
$\varepsilon.$ Then it could be shown that there exists a constant $u_{0},$
independent of $\varepsilon,$ such that, if $\left\Vert u_{\varepsilon
}\right\Vert _{L^{\infty}\left(  B_{\varepsilon}\right)  }\leq u_{0}$ then
there exists a unique solution $c_{\varepsilon}$ of problem (\ref{cd1}%
)$-$(\ref{cd2}). It satisfies an a priori estimate $\left\Vert c_{\varepsilon
}\right\Vert _{H^{1}\left(  B_{\varepsilon}\right)  }\leq c\left(  \left\Vert
g\right\Vert _{L^{2}\left(  B_{\varepsilon}\right)  }+\left\Vert G\right\Vert
_{H^{1}\left(  B_{\varepsilon}\right)  }\right)  ,$ where $c$ is independent
of $\varepsilon$ (see \cite{Lad}).

\subsection{The asymptotic expansion in a channel\label{sect1}}

Consider the equation
\begin{equation}
-\operatorname{div}\left(  \varkappa\nabla c_{\varepsilon r}\right)
+\mathbf{V}\left(  \frac{x_{2}}{\varepsilon}\right)  \cdot\nabla
c_{\varepsilon r}=g\left(  x_{1}\right)  , \label{cd1}%
\end{equation}
in the infinite channel $G_{\varepsilon}=\left(  -\infty,+\infty\right)
\times\left(  -\frac{\varepsilon}{2},\frac{\varepsilon}{2}\right)  ,$where
$\mathbf{V}\left(  \xi_{2}\right)  =\left(
\begin{array}
[c]{c}%
V_{p}\left(  \xi_{2}\right) \\
0
\end{array}
\right)  ,$ $V_{p}\left(  \xi_{2}\right)  =-\left(  \xi_{2}^{2}-\frac{1}%
{4}\right)  ;$ we consider the boundary conditions%
\begin{equation}
\varkappa\dfrac{\partial c_{\varepsilon r}}{\partial n}=\varepsilon\beta
c_{\varepsilon r}\text{ \ if \ }x_{2}=\pm\frac{\varepsilon}{2}. \label{cd2}%
\end{equation}
We will construct a function%
\begin{equation}
c_{\varepsilon r}^{(k)}\left(  x_{1},\xi_{2}\right)  =\sum_{j=0}%
^{k}\varepsilon^{j}c_{j}\left(  x_{1},\xi_{2}\right)  \label{cd3}%
\end{equation}
such that it satisfies (\ref{cd1}), (\ref{cd2}) up to the terms of order
$\varepsilon^{k}$ if $g\in C^{k+2}\left(  \mathbb{R}\right)  $. Substituting
(\ref{cd3}) into (\ref{cd1}), (\ref{cd2}) we get%
\[
\sum\limits_{j=0}^{k}\varepsilon^{j-2}\left(  \varkappa\frac{\partial^{2}%
c_{j}\left(  x_{1},\xi_{2}\right)  }{\partial\xi_{2}^{2}}+\varkappa
\frac{\partial^{2}c_{j-2}}{\partial x_{1}^{2}}-V_{p}\left(  \xi_{2}\right)
\dfrac{\partial c_{j-2}}{\partial x_{1}}\right)  =g\left(  x_{1}\right)
+R_{\varepsilon}^{\left(  k\right)  }%
\]
where $R_{\varepsilon}^{\left(  k\right)  }$ is a discrepancy%
\[
R_{\varepsilon}^{\left(  k\right)  }=\sum_{j=k-1,k}\varepsilon^{j}\left(
\varkappa\frac{\partial^{2}c_{j}}{\partial x_{1}^{2}}-V_{p}\left(  \xi
_{2}\right)  \dfrac{\partial c_{j}}{\partial x_{1}}\right)  ,
\]
and if $\xi_{2}=\pm\frac{1}{2},$
\[
\pm\sum\limits_{j=0}^{k}\varepsilon^{j-1}\varkappa\frac{\partial c_{j}%
}{\partial\xi_{2}}\left(  x_{1},\xi_{2}\right)  |_{\xi_{2}=\pm\frac{1}{2}%
}=\sum\limits_{j=0}^{k-1}\varepsilon^{j-1}\beta~c_{j-1}\left(  x_{1},\xi
_{2}\right)  +\varepsilon^{k}\beta~c_{k}\left(  x_{1},\xi_{2}\right)  .
\]
Equating the terms of the same power of $\varepsilon,$ we get%
\[
\varkappa\frac{\partial^{2}c_{j}}{\partial\xi_{2}^{2}}=-\varkappa
\dfrac{\partial^{2}c_{j-2}}{\partial x_{1}^{2}}+V_{p}\left(  \xi_{2}\right)
\dfrac{\partial c_{j-2}}{\partial x_{1}}+g\left(  x_{1}\right)  \delta_{j2}%
\]
and for $\xi_{2}=\pm\frac{1}{2}$%
\[
\pm\varkappa\frac{\partial c_{j}}{\partial\xi_{2}}|_{\xi_{2}=\pm\frac{1}{2}%
}=\beta~c_{j-2}\left(  x_{1},\xi_{2}\right)  .
\]
The necessary and sufficient condition of existence of $c_{j}$: the condition
\[
\varkappa\int_{-\frac{1}{2}}^{\frac{1}{2}}\frac{\partial^{2}c_{j}}{\partial
\xi_{2}^{2}}d\xi_{2}=\varkappa\frac{\partial c_{j}}{\partial\xi_{2}}\left(
x_{1},\frac{1}{2}\right)  -\varkappa\frac{\partial c_{j}}{\partial\xi_{2}%
}\left(  x_{1},-\frac{1}{2}\right)
\]
implies%
\[
-\varkappa\dfrac{\partial^{2}\left\langle c_{j-2}\right\rangle }{\partial
x_{1}^{2}}+\left\langle V_{p}\left(  \xi_{2}\right)  \dfrac{\partial c_{j-2}%
}{\partial x_{1}}\right\rangle +g\left(  x_{1}\right)  \delta_{j1}%
=\beta\left(  c_{j-2}\left(  x_{1},\frac{1}{2}\right)  +c_{j-2}\left(
x_{1},-\frac{1}{2}\right)  \right)
\]
where $\left\langle \cdot\right\rangle =\int_{-\frac{1}{2}}^{\frac{1}{2}}\cdot
d\xi_{2}.$

Each function $c_{j}$ is sought as a sum $c_{j}\left(  x_{1},\xi_{2}\right)
=\overline{c}_{j}\left(  x_{1}\right)  +\widetilde{c}_{j}\left(  x_{1},\xi
_{2}\right)  ,$ $\left\langle \widetilde{c}_{j}\right\rangle =0$ and
$\left\langle c_{j}\right\rangle =\overline{c}_{j}.$ So for $\overline
{c}_{j-2}$ we get the equation:
\[
-\varkappa\dfrac{\partial^{2}\overline{c}_{j-2}}{\partial x_{1}^{2}%
}+\left\langle V_{p}\right\rangle \dfrac{\partial\overline{c}_{j-2}}{\partial
x_{1}}+\left\langle V_{p}\dfrac{\partial\widetilde{c}_{j-2}}{\partial x_{1}%
}\right\rangle +g\left(  x_{1}\right)  \delta_{j2}=2\beta~\overline{c}%
_{j-2}+\beta\left(  \widetilde{c}_{j-2}\left(  x_{1},\frac{1}{2}\right)
+\widetilde{c}_{j-2}\left(  x_{1},-\frac{1}{2}\right)  \right)  .
\]
So, we get an algorithm of the successive determination of $\overline{c}_{j},$
$\widetilde{c}_{j}:$%
\[
j=0:\left\{
\begin{array}
[c]{c}%
\varkappa\dfrac{\partial^{2}\widetilde{c}_{0}}{\partial\xi_{2}^{2}}%
=0,\ \ \xi_{2}\in\left(  -\frac{1}{2},\frac{1}{2}\right) \\
\varkappa\dfrac{\partial\widetilde{c}_{0}}{\partial\xi_{2}}|_{\xi_{2}=\pm
\frac{1}{2}}=0,\ \ \left\langle \widetilde{c}_{0}\right\rangle =0
\end{array}
\right.
\]
and so, $\widetilde{c}_{0}=0;$%
\begin{equation}
-\varkappa\dfrac{\partial^{2}\overline{c}_{0}}{\partial x_{1}^{2}%
}+\left\langle V_{p}\right\rangle \dfrac{\partial\overline{c}_{0}}{\partial
x_{1}}+g\left(  x_{1}\right)  =2\beta~\overline{c}_{0}; \label{1.4000}%
\end{equation}
$j>1;$ find $\widetilde{c}_{j}:$%
\begin{equation}
\left\{
\begin{array}
[c]{l}%
\varkappa\dfrac{\partial^{2}\widetilde{c}_{j}}{\partial\xi_{2}^{2}}%
=-\varkappa\dfrac{\partial^{2}c_{j-2}}{\partial x_{1}^{2}}+V_{p}\left(
\xi_{2}\right)  \dfrac{\partial c_{j-2}}{\partial x_{1}}+g\left(
x_{1}\right)  \delta_{j1},\ \ \xi_{2}\in\left(  -\frac{1}{2},\frac{1}%
{2}\right) \\
\pm\varkappa\dfrac{\partial\widetilde{c}_{j}}{\partial\xi_{2}}|_{\xi_{2}%
=\pm\frac{1}{2}}=\beta c_{j-2}\left(  x_{1},\xi_{2}\right)  ,\ \left\langle
\widetilde{c}_{j}\right\rangle =0.
\end{array}
\right.  \label{1.100}%
\end{equation}
Remark $\widetilde{c}_{1}=0,$ and find $\overline{c}_{j}:$%
\[
-\varkappa\dfrac{\partial^{2}\overline{c}_{j}}{\partial x_{1}^{2}%
}+\left\langle V_{p}\right\rangle \dfrac{\partial\overline{c}_{j}}{\partial
x_{1}}+\left\langle V_{p}\dfrac{\partial\widetilde{c}_{j}}{\partial x_{1}%
}\right\rangle =2\beta\overline{c}_{j}+\beta\left(  \widetilde{c}_{j}\left(
x_{1},\frac{1}{2}\right)  +\widetilde{c}_{j}\left(  x_{1},-\frac{1}{2}\right)
\right)  .
\]

\begin{lemma}
\label{lem1}Each $\widetilde{c}_{j}\left(  x_{1},\xi_{2}\right)  $ is a
polynomial function in $\xi_{2}$ of degree $2j.$
\end{lemma}

\textbf{Proof.} By induction we prove that $\widetilde{c}_{j}$ is a polynomial
function of the degree $2j$. Indeed, $V_{p}$ is a quadratic function and so
the right hand side (\ref{1.100}) is a polynomial of order $2\left(
j-2\right)  +2.$ And so, after two integrations of (\ref{1.100}) we check that
$\widetilde{c}_{j}$ is a polynomial of order $2\left(  j-2\right)  +2+2=2j.$
$\blacksquare$

\bigskip

So we have constructed an asymptotic approximation (\ref{cd3}) which satisfies
equation (\ref{cd1}) up to the remainder $R^{\left(  k\right)  }$ and
conditions (\ref{cd2}) up to the remainder $\varepsilon^{k}\beta~c_{k}\left(
x_{1},\xi_{2}\right)  |_{\xi_{2}=\pm\frac{1}{2}}.$

Let us eliminate this remainder in the boundary conditions. To this end we
will add a corrector $\varepsilon^{k+1}\widehat{c}_{k+1}\left(  x_{1},\xi
_{2}\right)  $ such that
\[
\pm\varkappa\dfrac{\partial}{\partial x_{2}}\varepsilon^{k+1}\widehat{c}%
_{k+1}\left(  x_{1},\frac{x_{2}}{\varepsilon}\right)  =\varepsilon^{k+2}%
\beta\widehat{c}_{k+1}\left(  x_{1},\frac{x_{2}}{\varepsilon}\right)
+\varepsilon^{k}\beta c_{k}\left(  x_{1},\frac{x_{2}}{\varepsilon}\right)
\ \ \ \ \ \ \text{for }x_{2}=\pm\frac{\varepsilon}{2},
\]
i.e.
\begin{equation}
\pm\varkappa\dfrac{\partial}{\partial\xi_{2}}\varepsilon^{k}\widehat{c}%
_{k+1}\left(  x_{1},\xi_{2}\right)  =\varepsilon^{k+2}\beta\widehat{c}%
_{k+1}\left(  x_{1},\xi_{2}\right)  +\varepsilon^{k}\beta c_{k}\left(
x_{1},\xi_{2}\right)  \ \ \ \ \ \ \text{for }\xi_{2}=\pm\frac{1}{2}.
\label{1.200}%
\end{equation}
For example,%
\[
\widehat{c}_{k+1}=\varkappa^{-1}\left\{  \left(  \xi_{2}-\frac{1}{2}\right)
\left(  \xi_{2}+\frac{1}{2}\right)  ^{2}\beta c_{k}\left(  x_{1},\frac{1}%
{2}\right)  -\left(  \xi_{2}-\frac{1}{2}\right)  ^{2}\left(  \xi_{2}+\frac
{1}{2}\right)  \beta c_{k}\left(  x_{1},-\frac{1}{2}\right)  \right\}
\]
satisfies (\ref{1.200}). Then $c_{\varepsilon r}^{\left(  k\right)
}-\varepsilon^{k+1}\widehat{c}_{k+1}$ denoted by $\widehat{c}_{\varepsilon
r}^{\left(  k\right)  }$ satisfies the boundary conditions (\ref{cd2}) exactly
and the equation (\ref{cd1}) up to the remainder $R^{\left(  k\right)
}+\widehat{R}^{\left(  k\right)  },$ where
\[
\widehat{R}^{\left(  k\right)  }=\left\{  -div\left(  \varkappa\nabla
\widehat{c}_{k+1}\right)  -V_{p}\left(  \xi_{2}\right)  \frac{\partial
\widehat{c}_{k+1}}{\partial x_{1}}\right\}  \varepsilon^{k+1}.
\]

\subsection{Structural element "stenosis area".\label{sect2}}

Consider the convection-diffusion equation in the channel $G_{\varepsilon}$
with the modified coefficients $\varkappa$ and $V:$ they are replaced by the
functions
\[
K_{\varepsilon}\left(  x\right)  =\varkappa+\overline{K}\left(  \frac
{x}{\varepsilon}\right)  \text{ \ and \ }\mathbf{V}_{\varepsilon}\left(
x\right)  =\binom{V_{p}}{0}+\overline{V}\left(  \frac{x}{\varepsilon}\right)
\]
such that $K_{\varepsilon}\left(  x\right)  \geq\varkappa_{1}>0$ for all $x, $
$div\overline{V}\left(  \frac{x}{\varepsilon}\right)  =0$ and the support of
the functions $\overline{K}\left(  \xi\right)  $ and $\overline{V}\left(
\xi\right)  $ belongs to the ball $B\left(  O,2\right)  =\left\{  \xi
\in\mathbb{R}^{2}:|\xi|<2\right\}  .$

Consider the convection-diffusion equation%
\begin{equation}
-div\left(  K_{\varepsilon}\left(  x\right)  \nabla c_{\varepsilon}\right)
+\mathbf{V}_{\varepsilon}\left(  \frac{x}{\varepsilon}\right)  \cdot\nabla
c_{\varepsilon}=g_{\varepsilon}\left(  x_{1}\right)  ,\text{ \ in
}G_{\varepsilon} \label{1.300}%
\end{equation}
with the boundary conditions%
\begin{equation}
K_{\varepsilon}\left(  x\right)  \frac{\partial c_{\varepsilon}}{\partial
n}=\varepsilon\beta c_{\varepsilon},\ \ x_{2}=\pm\frac{\varepsilon}{2}.
\label{1.300bc}%
\end{equation}
Let us construct a function $c_{\varepsilon s}^{\left(  k\right)  }$
satisfying equation (\ref{1.300}) and conditions (\ref{1.300bc}) up to
remainder of order $\varepsilon^{k}.$ To this end we will consider the
asymptotic approximation $\widehat{c}_{\varepsilon r}^{\left(  k\right)  }$ of
the previous section completed by the special boundary layer corrector
$c_{\varepsilon sbl}^{\left(  k\right)  },$ having the form%
\[
c_{\varepsilon sbl}^{\left(  k\right)  }=\sum_{j=0}^{k}\varepsilon^{j}%
U_{j}\left(  \frac{x}{\varepsilon}\right)
\]
such that
\[
\left\vert U_{j}\left(  \xi\right)  \right\vert \leq c_{1}e^{-c_{2}\left\vert
\xi_{1}\right\vert }.
\]
Denote%
\[
L_{\varepsilon}=div\left(  K_{\varepsilon}\nabla\right)  -\mathbf{V}%
_{\varepsilon}\cdot\nabla.
\]
Define%
\[
\rho\left(  t\right)  =\left\{
\begin{array}
[c]{c}%
1\ \ \text{for}\ \left\vert t\right\vert >2;\\
0\ \ \text{for\ }\left\vert t\right\vert <1;
\end{array}
\right.
\]
such that $\rho$ and $\rho^{\prime}$ are bounded. Consider the following
asymptotic approximation%
\[
c_{\varepsilon s}^{(k)}=\widehat{c}_{\varepsilon r}^{\left(  k\right)  }%
\rho\left(  \frac{x_{1}}{\varepsilon}\right)  +\sum_{j=0}^{k}\varepsilon
^{j}U_{j}\left(  \frac{x}{\varepsilon}\right)  .
\]
Then%
\begin{align*}
L_{\varepsilon}c_{\varepsilon s}^{(k)}-g\left(  x_{1}\right)   &
=-L_{\varepsilon}\left(  \widehat{c}_{\varepsilon r}^{\left(  k\right)
}\left(  1-\rho\left(  \frac{x_{1}}{\varepsilon}\right)  \right)  \right)
-g+\left(  1-\rho\left(  \frac{x_{1}}{\varepsilon}\right)  \right)  \left(
L_{\varepsilon}\widehat{c}_{\varepsilon r}^{\left(  k\right)  }-g\right)  +\\
&  +L_{\varepsilon}\left(  \sum_{j=0}^{k}\varepsilon^{j}U_{j}\left(  \frac
{x}{\varepsilon}\right)  \right)  +\overline{R}_{k},
\end{align*}
where $\overline{R}_{k}=\rho\left(  \dfrac{x_{1}}{\varepsilon}\right)  \left(
L_{\varepsilon}\widehat{c}_{\varepsilon r}^{\left(  k\right)  }-g\right)  $
and $\left\vert \overline{R}_{k}\right\vert \leq\varepsilon^{k}const.$

Expand $\widehat{c}_{\varepsilon r}^{\left(  k\right)  }$ in powers of
$\varepsilon\xi_{1}:$%
\[
\widehat{c}_{\varepsilon r}^{\left(  k\right)  }=\sum_{j=0}^{k}\varepsilon
^{j}\sum_{l=0}^{k}\frac{\varepsilon^{l}}{l!}\xi_{1}^{l}\dfrac{\partial
^{l}c_{j}}{\partial x_{1}^{l}}\left(  0,\xi_{2}\right)  +\overline
{\overline{R}}_{k}=\sum_{r=0}^{k}\varepsilon^{r}\sum_{l=0}^{r}\frac{\xi
_{1}^{l}}{l!}\dfrac{\partial^{l}c_{r-l}}{\partial x_{1}^{l}}\left(  0,\xi
_{2}\right)  +\overline{\overline{\overline{R}}}_{k}.
\]
Let us denote $K\left(  \xi\right)  $ the function $\varkappa+\overline
{K}\left(  \xi\right)  $ and $\mathbf{V}\left(  \xi\right)  $ the function
$\binom{V_{p}}{0}+\overline{V}\left(  \xi\right)  .$ We have%
\begin{align*}
L_{\varepsilon}\left(  \widehat{c}_{\varepsilon r}^{\left(  k\right)  }\left(
1-\rho\left(  \frac{x_{1}}{\varepsilon}\right)  \right)  \right)   &
=\sum_{j=0}^{k}\varepsilon^{j-2}div_{\xi}\left(  K\left(  \xi\right)
\nabla_{\xi}\left(  \sum_{l=0}^{j}\frac{\xi_{1}^{l}}{l!}\dfrac{\partial
^{l}c_{j-l}}{\partial x_{1}^{l}}\left(  0,\xi_{2}\right)  \times\left(
1-\rho\left(  \xi_{1}\right)  \right)  \right)  \right)  -\\
&  -\sum_{j=0}^{k}\varepsilon^{j-1}\mathbf{V}\left(  \xi\right)  \cdot
\nabla_{\xi}\left(  \sum_{l=0}^{j}\frac{\xi_{1}^{l}}{l!}\dfrac{\partial
^{l}c_{j-l}}{\partial x_{1}^{l}}\left(  0,\xi_{2}\right)  \times\left(
1-\rho\left(  \xi_{1}\right)  \right)  \right)  .
\end{align*}
So, for $U_{j}$ we get equation%
\[
div_{\xi}\left(  K\left(  \xi\right)  \nabla_{\xi}U_{j}\left(  \xi\right)
\right)  -\mathbf{V}\left(  \xi\right)  \cdot\nabla_{\xi}U_{j-1}\left(
\xi\right)  +T_{j}\left(  \xi\right)  =0,\ \ \xi_{1}\in\mathbb{R},\ \xi_{2}%
\in\left(  -\frac{1}{2},\frac{1}{2}\right)  ,
\]
where
\begin{align*}
T_{j}\left(  \xi\right)   &  =div_{\xi}\left(  K\left(  \xi\right)
\nabla_{\xi}\left(  \sum_{l=0}^{j}\frac{\xi_{1}^{l}}{l!}\dfrac{\partial
^{l}c_{j-l}}{\partial x_{1}^{l}}\left(  0,\xi_{2}\right)  \times\left(
1-\rho\left(  \xi_{1}\right)  \right)  \right)  \right)  -\\
&  -\mathbf{V}\left(  \xi\right)  \cdot\nabla_{\xi}\left(  \sum_{l=0}%
^{j-1}\frac{\xi_{1}^{l}}{l!}\dfrac{\partial^{l}c_{j-1-l}}{\partial x_{1}^{l}%
}\left(  0,\xi_{2}\right)  \times\left(  1-\rho\left(  \xi_{1}\right)
\right)  \right)  .
\end{align*}
For the boundary conditions, in the same way:%
\[
\pm\sum_{j=0}^{k}\varepsilon^{j-1}K\left(  \xi\right)  \dfrac{\partial U_{j}%
}{\partial\xi_{2}}\pm\sum_{j=0}^{k}\varepsilon^{j-1}K\left(  \xi\right)
\dfrac{\partial c_{j}}{\partial\xi_{2}}\left(  \rho\left(  \xi_{1}\right)
-1\right)  =\varepsilon\beta\left(  \sum_{j=0}^{k}\varepsilon^{j}U_{j}%
+\sum_{j=0}^{k}\varepsilon^{j}c_{j}\left(  \rho\left(  \xi_{1}\right)
-1\right)  \right)  ,\ \ \xi_{2}=\pm\frac{1}{2}%
\]
and so%
\[
\pm K\left(  \xi\right)  \dfrac{\partial U_{j}}{\partial\xi_{2}}=S_{j}\left(
\xi\right)  ,\ \ \xi_{2}=\pm\frac{1}{2}%
\]
where%
\[
S_{j}\left(  \xi\right)  =-\left\{  \pm\sum_{l=0}^{j}K\left(  \xi\right)
\frac{\xi_{1}^{l}}{l!}\dfrac{\partial^{l+1}c_{j-l}}{\partial\xi_{2}\partial
x_{1}^{l}}\left(  0,\xi_{2}\right)  \left(  1-\rho\left(  \xi_{1}\right)
\right)  -\beta\sum_{l=0}^{j-2}\frac{\xi_{1}^{l}}{l!}\dfrac{\partial^{l}c_{j}%
}{\partial x_{1}^{l}}\left(  0,\xi_{2}\right)  \left(  1-\rho\left(  \xi
_{1}\right)  \right)  -\beta U_{j-2}\right\}  .
\]
Necessary and sufficient condition of existence of a bounded solution $U_{j}:
$%
\[
\left\{  T_{j}\left(  \xi\right)  -\mathbf{V}\left(  \xi\right)  \cdot
\nabla_{\xi}U_{j-1}\right\}  =\left\{  S_{j}\right\}  _{+}+\left\{
S_{j}\right\}  _{-}%
\]
where%
\[
\left\{  \cdot\right\}  =\int_{\left(  -\infty,0\right)  \times\left(
-\frac{1}{2},\frac{1}{2}\right)  }\cdot d\xi+\int_{\left(  0,+\infty\right)
\times\left(  -\frac{1}{2},\frac{1}{2}\right)  }\cdot d\xi
\]
and
\[
\left\{  \cdot\right\}  _{\pm}=\int_{\left(  -\infty,0\right)  }\cdot
|_{\xi_{2}=\pm\frac{1}{2}}d\xi_{1}+\int_{\left(  0,+\infty\right)  }%
\cdot|_{\xi_{2}=\pm\frac{1}{2}}d\xi_{1}.
\]
This condition gives one interface condition for $\dfrac{\partial\overline
{c}_{j-1}}{\partial x_{1}}\left(  0^{+}\right)  $ and $\dfrac{\partial
\overline{c}_{j-1}}{\partial x_{1}}\left(  0^{-}\right)  :$%
\[%
\begin{array}
[c]{l}%
\left\{  T_{j}\right\}  =%
%TCIMACRO{\dint _{-\frac{1}{2}}^{\frac{1}{2}}}%
%BeginExpansion
{\displaystyle\int_{-\frac{1}{2}}^{\frac{1}{2}}}
%EndExpansion
K\left(  \xi\right)  d\xi_{2}\dfrac{\partial\overline{c}_{j-1}}{\partial
x_{1}}\left(  0^{-}\right)  -%
%TCIMACRO{\dint _{-\frac{1}{2}}^{\frac{1}{2}}}%
%BeginExpansion
{\displaystyle\int_{-\frac{1}{2}}^{\frac{1}{2}}}
%EndExpansion
K\left(  \xi\right)  d\xi_{2}\dfrac{\partial\overline{c}_{j-1}}{\partial
x_{1}}\left(  0^{+}\right)  +\\
\ \ \ \ \ \ \ \ \ +%
%TCIMACRO{\dint _{-\frac{1}{2}}^{\frac{1}{2}}}%
%BeginExpansion
{\displaystyle\int_{-\frac{1}{2}}^{\frac{1}{2}}}
%EndExpansion
K\left(  \xi\right)  \dfrac{\partial\widetilde{c}_{j-1}}{\partial x_{1}%
}\left(  0^{-},\xi_{2}\right)  d\xi_{2}-%
%TCIMACRO{\dint _{-\frac{1}{2}}^{\frac{1}{2}}}%
%BeginExpansion
{\displaystyle\int_{-\frac{1}{2}}^{\frac{1}{2}}}
%EndExpansion
K\left(  \xi\right)  \dfrac{\partial\widetilde{c}_{j-1}}{\partial x_{1}%
}\left(  0^{+},\xi_{2}\right)  d\xi_{2}+\\
\ \ \ \ \ \ \ \ \ +\ \text{all other terms of }T_{j}\text{ except }l=1\text{
for }div.
\end{array}
\]

So, if we pose $Q_{j}=\int_{-\frac{1}{2}}^{\frac{1}{2}}K\left(  \xi\right)
\dfrac{\partial\widetilde{c}_{j-1}}{\partial x_{1}}\left(  0^{-},\xi
_{2}\right)  d\xi_{2}-\int_{-\frac{1}{2}}^{\frac{1}{2}}K\left(  \xi\right)
\dfrac{\partial\widetilde{c}_{j-1}}{\partial x_{1}}\left(  0^{+},\xi
_{2}\right)  d\xi_{2},$ we get condition:%
\begin{equation}
-\int_{-\frac{1}{2}}^{\frac{1}{2}}K\left(  0,\xi_{2}\right)  d\xi_{2}\left[
\dfrac{\partial\overline{c}_{j-1}}{\partial x_{1}}\right]  =g_{j}
\label{1.400}%
\end{equation}
where
\[
g_{j}=-\left\{  T_{j}\left(  \xi\right)  -\mathbf{V}\left(  \xi\right)
\cdot\nabla_{\xi}U_{j-1}\right\}  -Q_{j}+\left\{  S_{j}\right\}  _{+}+\left\{
S_{j}\right\}  _{-},
\]
where in $T_{j}\left(  \xi\right)  $ there are not the term corresponding to
$l=1$ for $div$.

Now we solve the problem on $U_{j}$ exponentially stabilizing to some
constants. We choose this constant equal to $0$ at $-\infty.$ At $+\infty$ we
have $\widetilde{U}_{j}\rightarrow\widetilde{q}_{j}.$ To make $U_{j}%
\rightarrow0$ as $\xi_{1}\rightarrow+\infty$, we subtract this constant for
all $\xi_{1}>0;$ we set%
\[
U_{j}=\left\{
\begin{array}
[c]{l}%
\widetilde{U}_{j}\ \ \text{for }\xi_{1}<0\\
\widetilde{U}_{j}-\widetilde{q}_{j}\ \ \text{for }\xi_{1}>0.
\end{array}
\right.
\]
This function $U_{j}$ is exponentially decaying at $\infty.$ Then we put the
compensating condition for $\overline{c}_{j}$ at $x_{1}=0:\left[  \overline
{c}_{j}\right]  =\widetilde{q}_{j}.$

Another interface condition for $\overline{c}_{j}$ is (\ref{1.400}).

\subsection{The boundary layer in a "bifurcation area".\label{sect3}}

Consider the convection-diffusion equation in a one bundle structure with the
common point $O$ (i.e. the origin of the coordinate system). We will construct
an asymptotic expansion in the neighborhood of the point $O$ (that will be
good in all the channels $B_{i}^{\varepsilon})$. First we construct a regular
expansion as in section \ref{sect1} for every branch $B_{i}^{\varepsilon}$
making the local change of variables by rotation in such a way that new
$x_{1}$ axis (denoted $x_{1}^{e_{i}})$ contains the segment $e_{i}.$ Let us
denote the approximation $\widehat{c}_{\varepsilon r}^{\left(  k\right)  }$
for $e_{i}$ as%
\[
\widehat{c}_{\varepsilon r}^{\left(  k\right)  e_{i}}\left(  x\right)
=\sum_{j=0}^{k}\varepsilon^{j}c_{j}\left(  x_{1}^{e_{i}},\frac{x_{2}^{e_{i}}%
}{\theta_{i}\varepsilon}\right)  -\varepsilon^{k+1}\widehat{c}_{k+1}\left(
x_{1}^{e_{i}},\frac{x_{2}^{e_{i}}}{\theta_{i}\varepsilon}\right)  .
\]
Here $\left(  x_{1}^{e_{i}},x_{2}^{e_{i}}\right)  $ are new local variables
related to the segment $e_{i},$ and $\theta_{i}\varepsilon$ is the thickness
of the channel $B_{i}^{\varepsilon}.$ Then, as in section \ref{sect2}, we
construct the boundary layer corrector%
\[
c_{\varepsilon bl}^{(k)}=\sum_{j=0}^{k}\varepsilon^{j}U_{j}\left(  \frac
{x}{\varepsilon}\right)
\]
such that
\[
\left\vert U_{j}\left(  \xi\right)  \right\vert \leq\overline{c}%
_{1}e^{-\overline{c}_{2}\left\vert \xi\right\vert },\ \ \overline{c}%
_{1},\overline{c}_{2}>0.
\]
Consider an asymptotic approximation%
\begin{equation}
c_{\varepsilon b}^{\left(  k\right)  }=\sum_{i=1}^{n}\widehat{c}_{\varepsilon
r}^{\left(  k\right)  }\left(  x_{1}^{e_{i}},\frac{x_{2}^{e_{i}}}{\theta
_{i}\varepsilon}\right)  \chi_{i}\left(  x\right)  \rho\left(  \frac{\alpha
x_{1}^{e_{i}}}{\varepsilon}\right)  +\sum_{j=0}^{k}\varepsilon^{j}U_{j}\left(
\frac{x}{\varepsilon}\right)  , \label{1.5000}%
\end{equation}
where $\chi_{i}$ is a characteristic function of $B_{i}^{\varepsilon}$
($\chi_{i}\left(  x\right)  =1$ if $x\in B_{i}^{\varepsilon},$ and $0$ if
not); $\alpha$ is a number such that each point of the circle of a radius
$\alpha\varepsilon$ with center $O$ belongs to at most one rectangle
$B_{i}^{\varepsilon}$ ($i=1,...,n)$ and it does not intersect the domain
$\gamma_{O}^{\varepsilon}.$ It means that the variation of the function $\rho$
from 1 to 0 occurs out of this circle.%
\[%
%TCIMACRO{\FRAME{itbpFU}{1.9277in}{1.8533in}{0in}{\Qcb{Fig. 3}}{}%
%{fig2.eps}{\special{ language "Scientific Word";  type "GRAPHIC";
%width 1.9277in;  height 1.8533in;  depth 0in;  original-width 2.5685in;
%original-height 2.469in;  cropleft "0";  croptop "1";  cropright "1";
%cropbottom "0";  filename 'fig2.eps';file-properties "XNPEU";}} }%
%BeginExpansion
{\parbox[b]{1.9277in}{\begin{center}
\includegraphics[
height=1.8533in,
width=1.9277in
]%
{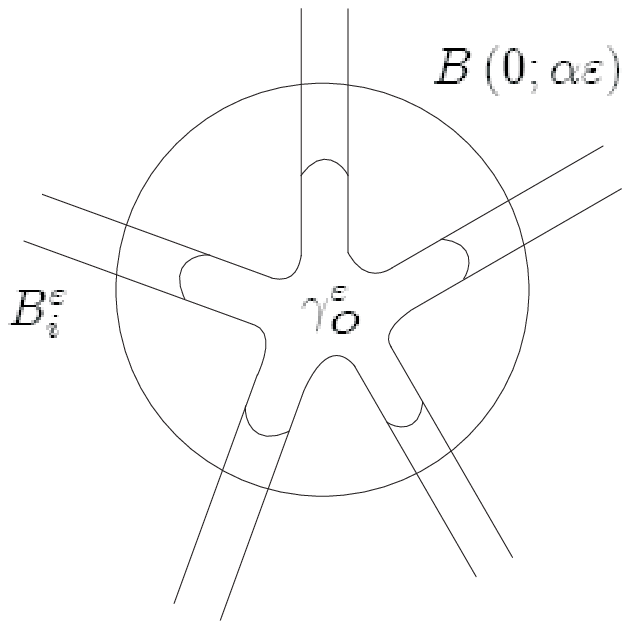}%
\\
Fig. 3
\end{center}}}
%EndExpansion
\]

Assume that the right hand side is equal to zero in some neighborhood of $O.$
Then for every $i=1,...,n$ for $x\in B_{i}^{\varepsilon}$%
\begin{align*}
L_{\varepsilon}c_{\varepsilon b}^{(k)e_{i}}-g\left(  x_{1}^{e_{i}}\right)   &
=-L_{\varepsilon}\left(  \widehat{c}_{\varepsilon r}^{\left(  k\right)  e_{i}%
}\left(  1-\rho\left(  \frac{\alpha x_{1}^{e_{i}}}{\varepsilon}\right)
\right)  \right)  +\left(  1-\rho\left(  \frac{\alpha x_{1}^{e_{i}}%
}{\varepsilon}\right)  \right)  \left(  L_{\varepsilon}\widehat{c}%
_{\varepsilon r}^{\left(  k\right)  e_{i}}-g\right)  +\\
&  +L_{\varepsilon}\left(  \sum_{j=0}^{k}\varepsilon^{j}U_{j}\left(  \frac
{x}{\varepsilon}\right)  \right)  +\overline{R}_{k},
\end{align*}
where $\overline{R}_{k}=\rho\left(  \dfrac{\alpha x_{1}^{e_{i}}}{\varepsilon
}\right)  \left(  L_{\varepsilon}\widehat{c}_{\varepsilon r}^{\left(
k\right)  e_{i}}-g\right)  .$ Inside the circle $B\left(  0;\alpha
\varepsilon\right)  $ we get: $g=0,$ $\rho=0$ and so
\[
L_{\varepsilon}c_{\varepsilon b}^{(k)e_{i}}-g=L_{\varepsilon}\left(
\sum_{j=0}^{k}\varepsilon^{j}U_{j}\left(  \frac{x}{\varepsilon}\right)
\right)  .
\]
Here the terms $\left(  1-\rho\left(  \dfrac{\alpha x_{1}^{e_{i}}}%
{\varepsilon}\right)  \right)  \left(  L_{\varepsilon}\widehat{c}_{\varepsilon
r}^{\left(  k\right)  e_{i}}-g\right)  $ and $\rho\left(  \dfrac{\alpha
x_{1}^{e_{i}}}{\varepsilon}\right)  \left(  L_{\varepsilon}\widehat
{c}_{\varepsilon r}^{\left(  k\right)  e_{i}}-g\right)  $ are bounded by
$\varepsilon^{k}const.$ So we will define $U_{j}$ in such a way that out of
the circle $B\left(  0;\alpha\varepsilon\right)  $ we have%
\[
L_{\varepsilon}\left(  \sum_{j=0}^{k}\varepsilon^{j}U_{j}\left(  \frac
{x}{\varepsilon}\right)  \right)  -\sum_{i=1}^{n}L_{\varepsilon}\left(
\widehat{c}_{\varepsilon r}^{\left(  k\right)  e_{i}}\chi_{i}\left(
1-\rho\left(  \frac{\alpha x_{1}^{e_{i}}}{\varepsilon}\right)  \right)
\right)
\]
to be equal to zero up to the terms of order $\varepsilon^{k}$ and inside the
circle,%
\[
L_{\varepsilon}\left(  \sum_{j=0}^{k}\varepsilon^{j}U_{j}\left(  \frac
{x}{\varepsilon}\right)  \right)  =0.
\]
Expand each $\widehat{c}_{\varepsilon r}^{\left(  k\right)  e_{i}}$ in powers
of $\varepsilon\xi_{1}^{e_{i}}:$%
\begin{align*}
\widehat{c}_{\varepsilon r}^{\left(  k\right)  e_{i}}  &  =\sum_{j=0}%
^{k}\varepsilon^{j}\sum_{l=0}^{j}\frac{\varepsilon^{l}}{l!}\left(  \xi
_{1}^{e_{i}}\right)  ^{l}\dfrac{\partial^{l}c_{j}^{e_{i}}}{\partial\left(
x_{1}^{e_{i}}\right)  ^{l}}\left(  0,\xi_{2}^{e_{i}}\right)  +\overline
{\overline{R}}_{k}=\\
&  =\sum_{r=0}^{k}\varepsilon^{r}\sum_{l=0}^{r}\frac{\xi_{1}^{l}}{l!}%
\dfrac{\partial^{l}c_{r-l}^{e_{i}}}{\partial\left(  x_{1}^{e_{i}}\right)
^{l}}\left(  0,\xi_{2}^{e_{i}}\right)  +\overline{\overline{\overline{R}}}%
_{k},\ \ \xi_{2}^{e_{i}}=\frac{x_{2}^{e_{i}}}{\varepsilon}.
\end{align*}
Consider domain $\Omega_{0}=\gamma_{0}\cup%
%TCIMACRO{\tbigcup \limits_{i=1}^{n}}%
%BeginExpansion
{\textstyle\bigcup\limits_{i=1}^{n}}
%EndExpansion
\Omega_{i},$ where $\Omega_{i}$ is the half-strip $\left\{  \left(  \xi
_{1},\xi_{2}\right)  \in\mathbb{R}^{2}:\xi_{1}^{e_{i}}>0,\left\vert \xi
_{2}^{e_{i}}\right\vert <\dfrac{\theta_{i}}{2}\right\}  .$ So for $U_{j}$ we
get equation:%

\[
\varkappa\Delta U_{j}\left(  \xi\right)  -\mathbf{V}\left(  \xi\right)
\cdot\nabla_{\xi}U_{j-1}\left(  \xi\right)  +T_{j}\left(  \xi\right)
=0,\ \ \xi\in\Omega_{0},
\]
where
\[
T_{j}\left(  \xi\right)  =\left\{
\begin{array}
[c]{l}%
T_{j}^{e_{i}}\left(  \xi\right)  \ \ \ \text{if \ }\xi\in\Omega_{i}\text{ out
of the circle }B\left(  0,\alpha\right)  ,\\
0\text{ \ \ inside the circle }B\left(  0,\alpha\right)  ,
\end{array}
\right.
\]%
\begin{align*}
T_{j}^{e_{i}}\left(  \xi\right)   &  =\varkappa\Delta_{\xi}\left(  \sum
_{l=0}^{j}\frac{\left(  \xi_{1}^{e_{i}}\right)  ^{l}}{l!}\dfrac{\partial
^{l}c_{j-l}}{\partial\left(  x_{1}^{e_{i}}\right)  ^{l}}\left(  0,\xi
_{2}^{e_{i}}\right)  \times\left(  1-\rho\left(  \alpha\xi_{1}^{e_{i}}\right)
\right)  \right)  -\\
&  -\mathbf{V}\left(  \xi\right)  \cdot\nabla_{\xi}\left(  \sum_{l=0}^{j}%
\frac{\left(  \xi_{1}^{e_{i}}\right)  ^{l}}{l!}\dfrac{\partial^{l}c_{j-1-l}%
}{\partial\left(  x_{1}^{e_{i}}\right)  ^{l}}\left(  0,\xi_{2}^{e_{i}}\right)
\right)  .
\end{align*}
For the boundary conditions on $\partial\Omega_{0}$:%
\[%
\begin{array}
[c]{c}%
\pm%
%TCIMACRO{\dsum \limits_{j=0}^{k}}%
%BeginExpansion
{\displaystyle\sum\limits_{j=0}^{k}}
%EndExpansion
\varepsilon^{j-1}\varkappa\dfrac{\partial U_{j}}{\partial\xi_{2}^{e_{i}}}\pm%
%TCIMACRO{\dsum \limits_{j=0}^{k}}%
%BeginExpansion
{\displaystyle\sum\limits_{j=0}^{k}}
%EndExpansion
\varepsilon^{j-1}\varkappa\dfrac{\partial c_{j}^{e_{i}}}{\partial\xi
_{2}^{e_{i}}}\left(  x_{1}^{e_{i}},\xi_{2}^{e_{i}}\right)  \left(  \rho\left(
\alpha\xi_{1}^{e_{i}}\right)  -1\right)
=\ \ \ \ \ \ \ \ \ \ \ \ \ \ \ \ \ \ \ \ \ \ \ \ \ \ \ \ \ \ \ \ \ \ \ \ \ \ \ \ \ \ \ \ \ \ \ \ \ \ \ \ \ \ \ \ \ \ \ \ \ \ \ \ \ \ \ \ \ \ \ \ \ \ \ \\
=\varepsilon\beta%
%TCIMACRO{\dsum \limits_{j=0}^{k}}%
%BeginExpansion
{\displaystyle\sum\limits_{j=0}^{k}}
%EndExpansion
\varepsilon^{j}U_{j}+%
%TCIMACRO{\dsum \limits_{j=0}^{k}}%
%BeginExpansion
{\displaystyle\sum\limits_{j=0}^{k}}
%EndExpansion
\varepsilon^{j}c_{j}^{e_{i}}\left(  \rho\left(  \alpha\xi_{1}^{e_{i}}\right)
-1\right)  ,\ \text{out of the circle }B\left(  0,\alpha\right)
\end{array}
\]
and%
\[
\pm\sum_{j=0}^{k}\varepsilon^{j-1}\varkappa\dfrac{\partial U_{j}}{\partial
n_{\xi}}=\varepsilon\beta\sum_{j=0}^{k}\varepsilon^{j}U_{j},\text{ inside the
circle }B\left(  0,\alpha\right)
\]
up to the terms of order $\varepsilon^{k}.$ This gives%
\[
\varkappa\dfrac{\partial U_{j}}{\partial n_{\xi}}=S_{j}\left(  \xi\right)
,\ \ \xi\in\partial\Omega_{0},
\]
where%
\[
S_{j}\left(  \xi\right)  =\left\{
\begin{array}
[c]{l}%
\begin{array}
[c]{c}%
-\left\{  \pm%
%TCIMACRO{\dsum \limits_{l=0}^{j}}%
%BeginExpansion
{\displaystyle\sum\limits_{l=0}^{j}}
%EndExpansion
\varkappa\dfrac{\left(  \xi_{1}^{e_{i}}\right)  ^{l}}{l!}\dfrac{\partial
^{l+1}c_{j-l}^{e_{i}}}{\partial\xi_{2}^{e_{i}}\partial\left(  x_{1}^{e_{i}%
}\right)  ^{l}}\left(  0,\xi_{2}^{e_{i}}\right)  \left(  1-\rho\left(
\alpha\xi_{1}^{e_{i}}\right)  \right)  -\right.
\ \ \ \ \ \ \ \ \ \ \ \ \ \ \ \ \ \ \ \ \ \ \ \ \ \ \ \ \ \ \ \ \ \ \ \ \ \ \ \ \ \ \ \ \ \ \ \ \ \ \ \ \ \ \ \ \\
\left.  -\beta%
%TCIMACRO{\dsum \limits_{l=0}^{j-2}}%
%BeginExpansion
{\displaystyle\sum\limits_{l=0}^{j-2}}
%EndExpansion
\dfrac{\left(  \xi_{1}^{e_{i}}\right)  ^{l}}{l!}\dfrac{\partial^{l}%
c_{j}^{e_{i}}}{\partial\left(  x_{1}^{e_{i}}\right)  ^{l}}\left(  0,\xi
_{2}^{e_{i}}\right)  \left(  1-\rho\left(  \alpha\xi_{1}^{e_{i}}\right)
\right)  -\beta U_{j-2},\right\}  \ \text{out of the circle }B\left(
0,\alpha\right)  \text{ on }\partial\Omega_{i},
\end{array}
\\
\\
+\beta U_{j-2},\text{ \ inside the circle }B\left(  0,\alpha\right)  \text{ on
}\partial\Omega_{0}%
\end{array}
\right.  .
\]
Necessary and sufficient condition of existence of a bounded solution $U_{j}:
$%
\[
\left\{  T_{j}\left(  \xi\right)  -\mathbf{V}\left(  \xi\right)  \cdot
\nabla_{\xi}U_{j-1}\right\}  _{\Omega_{0}}=\left\{  S_{j}\right\}
_{\partial\Omega_{0}},
\]
where%
\[
\left\{  \cdot\right\}  _{\Omega_{0}}=\int_{\Omega_{0}}\cdot d\xi\text{ \ and
\ }\left\{  \cdot\right\}  _{\partial\Omega_{0}}=\int_{\partial\Omega_{0}%
}\cdot ds.
\]
This condition gives one interface condition for $\dfrac{\partial\overline
{c}^{e_{i}}}{\partial x_{1}^{e_{i}}}\left(  0^{e_{i}}\right)  $, the limit
value of the derivative $\dfrac{\partial\overline{c}^{e_{i}}}{\partial
x_{1}^{e_{i}}}$ on $e_{i}$ in the origin of $e_{i}:$%
\[
\left\{  T_{j}\right\}  _{\Omega_{0}}=-\sum_{i=1}^{n}\varkappa\theta_{i}%
\dfrac{\partial\overline{c}_{j-1}^{e_{i}}}{\partial x_{1}^{e_{i}}}\left(
0^{e_{i}}\right)  -\sum_{i=1}^{n}\varkappa\theta_{i}\dfrac{\partial
\widetilde{c}_{j-1}^{e_{i}}}{\partial x_{1}^{e_{i}}}\left(  0^{e_{i}},\xi
_{2}^{e_{i}}\right)  +\text{all other terms of }T_{j}\text{ except }l=1\text{
for }\Delta_{\xi}.
\]
So, we get Kirchoff type condition:%
\[
-\varkappa\sum_{i=1}^{n}\theta_{i}\dfrac{\partial\overline{c}_{j-1}}{\partial
x_{1}^{e_{i}}}\left(  0^{e_{i}}\right)  =g_{j}%
\]
where $g_{j}$ depends on $\overline{c}_{0},...,\overline{c}_{j-2}.$

Now we solve the problem on $U_{j}$ in $\Omega_{0};$ it stabilizes
exponentially to some constants $\widetilde{q}_{ji}$ at every branch
$\Omega_{i}$. This solution $U_{j}$ corresponds to some values of
$\overline{c}_{j}^{e_{i}}\left(  0^{e_{i}}\right)  $ that enter in the
expression for $T_{j}^{e_{i}}$ (when $l=0$). If we change these values adding
$\widetilde{q}_{ji}$ then the solution $U_{j}$ will be transformed into
$U_{j}-\rho\left(  \alpha\xi_{1}^{e_{i}}\right)  \widetilde{q}_{ji}$ which
tend to zero as $\left\vert \xi\right\vert \rightarrow+\infty.$

It means that for every $j$ we have to solve first the problem for
$\overline{c}_{j},$ that is, an ordinary differential equation on every
$e_{i},$ Kirchoff condition for $\dfrac{\partial\overline{c}_{j-1}}{\partial
x_{1}^{e_{i}}}$ in $0$ and continuity condition for $\overline{c}_{j}$ in 0,
solve the problem for $U_{j}$ and then modify the values $\overline{c}%
_{j}\left(  0^{e_{i}}\right)  $ in such a way that $U_{j}\rightarrow0$ as
$\left\vert \xi\right\vert \rightarrow+\infty.$

\subsection{The entrance/exit element}

In this case the boundary layer is constructed in the following way. It
satisfies equation (\ref{cd1}), boundary conditions (\ref{cd2}) everywhere
except the part $\Gamma_{s},$ where $U_{j}=q_{s}\delta_{0j}.$ There is no more
necessary and sufficient condition for the existence of a bounded solution.
This bounded solution exponentially stabilizes to some constant. Then we
redefine the value $\overline{c}_{j}\left(  x_{bs}\right)  $ in such a way
that $U\rightarrow0$ as $\xi_{1}^{e_{i}}\rightarrow+\infty.$ In particular,
$\overline{c}_{0}\left(  x_{bs}\right)  =q_{s}$ and $U_{0}=\left(
1-\rho\right)  q_{s}.$

\subsection{Algorithm of assembling of an asymptotic solution}

Let us describe more precisely the algorithm of the assembling of the
different structural elements in the construction of the asymptotic expansion.

For $\overline{c}_{0}$ we get from the construction of the subsection
\ref{sect1} the 1D convection diffusion equation
\begin{align}
\varkappa\theta_{i}\dfrac{\partial^{2}\overline{c}_{0}^{e_{i}}}{\partial
\left(  x_{1}^{e_{i}}\right)  ^{2}}-\left\langle V_{p}\right\rangle
_{\theta_{i}}\dfrac{\partial\overline{c}_{0}^{e_{i}}}{\partial x_{1}^{e_{i}}%
}+\theta_{i}g\left(  x_{1}\right)   &  =2\beta\overline{c}_{0}^{e_{i}},\text{
\ for every }e_{i},\label{1.500}\\
\left\langle V_{p}\right\rangle _{\theta_{i}}  &  =\int_{-\frac{\theta_{i}}%
{2}}^{\frac{\theta_{i}}{2}}V_{p}\left(  \xi_{2}^{e_{i}}\right)  d\xi
_{2}^{e_{i}}.\nonumber
\end{align}
At every stenosis point $x_{s}$ we get the interface conditions%
\begin{equation}
\left[  \varkappa\dfrac{\partial\overline{c}_{0}^{e_{i}}}{\partial
x_{1}^{e_{i}}}\right]  =0\text{ \ and \ }\left[  \overline{c}_{0}^{e_{i}%
}\right]  =0. \label{1.6000}%
\end{equation}
At every bifurcation point $x_{b}$ that is an end point of segments
$e_{1},...,e_{n}$ we get%
\begin{equation}
\varkappa\sum_{i=1}^{n}\theta_{i}\dfrac{\partial\overline{c}_{0}^{e_{i}}%
}{\partial x_{1}^{e_{i}}}\left(  x_{b}\right)  =0 \label{1.6100}%
\end{equation}
and $\overline{c}_{0}^{e_{i}}\left(  x_{b}\right)  =\overline{c}_{0}^{e_{j}%
}\left(  x_{b}\right)  $ $\forall i,j\in\left\{  1,...,n\right\}  $ (continuity).

At every entrance/exit point $x_{t}$ we set $\overline{c}_{0}^{e_{i}}\left(
x_{t}\right)  =q_{t},$ $t=1,..,r,$ where $x_{t}$ is an end-point corresponding
to $\Gamma_{t}.$

Then for every $j$ we get analogous problem for $\overline{c}_{j}$ with some
right hand sides depending on previous approximations $\overline{c}%
_{0},...,\overline{c}_{j-1}.$ We solve the problem for $\widetilde{c}_{j}$ and
then for $U_{j}.$ We define the interface values $\overline{c}_{j}^{e_{i}%
}\left(  x_{b}\right)  $ for all segments $e_{i}$ having $x_{b}$ as an end
point in such a way that $U_{j}\rightarrow0$ for $\left\vert \xi\right\vert
\rightarrow+\infty.$

\begin{remark}
Every function $\overline{c}_{0}^{e_{i}}$ related to $e_{i}$ depends on the
local variable $x_{1}^{e_{i}}$ such that $x_{1}^{e_{i}}=0$ in one of the end
points of $e_{i}.$ Therefore we get two local variables $x_{1}^{e_{i}}$
"starting" from each of two ends of the segment $e_{i};$ therefore, if we
denote $x_{1}^{e_{i}+}$ and $x_{1}^{e_{i}-}$ these two different local
variables, we get the relation for these variables $x_{1}^{e_{i}+}=\left\vert
e_{i}\right\vert -x_{1}^{e_{i}-}.$ Let us give more details. Every segment
$e_{i}$ having two end points $O_{i_{1}}$ and $O_{i_{2}}$ can be associated to
two possible local coordinate systems: one of them (denoted by $Ox^{e_{i_{1}}%
}$) has its origin in $O_{i_{1}}$ and another in $O_{i_{2}}$ (denoted by
$Ox^{e_{i_{2}}}$). Therefore for every differential equation (\ref{1.500}) or
for every junction condition this system should be chosen and fixed. The
evident change variable relation for any function $f$ defined on the segment
$e_{i}$ is:%
\[
f\left(  x_{1}^{e_{i_{1}}}\right)  =f\left(  \left\vert e_{i}\right\vert
-x_{1}^{e_{i_{2}}}\right)  ,
\]
where $\left\vert e_{i}\right\vert $ is the length of the segment, in the left
side we use the first variable of the system $Ox^{e_{i_{1}}}$ and in the right
side $-$ the first variable of the system $Ox^{e_{i_{2}}}.$ Consequently, for
$x=O_{i_{1}},$ we have $\dfrac{\partial f}{\partial x_{1}^{e_{i_{1}}}}\left(
0\right)  =-\dfrac{\partial f}{\partial x_{1}^{e_{i_{2}}}}\left(  \left\vert
e_{i}\right\vert \right)  .$
\end{remark}

\subsection{Justification: draft.}

Substituting the asymptotic solution of a form (\ref{1.5000}) we satisfy the
convection-diffusion equation with a discrepancy $O\left(  \varepsilon
^{k}\right)  $ in $L^{2}$ norm. The boundary conditions are satisfied with the
same accuracy. Then applying the Poincar\'{e}-Friedrichs inequality for rod
structures (see \cite{GPbook}) and the a priori estimate derived from the
variational formulation, we get that%
\[
\left\Vert c_{\varepsilon}-c_{\varepsilon}^{(k)}\right\Vert _{H^{1}\left(
B_{\varepsilon}\right)  }\leq C\sqrt{\varepsilon}\varepsilon^{k},
\]
where $C$ does not depend on $\varepsilon.$

In particular the leading term $c_{0}$ satisfies the following estimate
\[
\left\Vert c_{\varepsilon}-c_{0}\right\Vert _{L^{2}\left(  B_{\varepsilon
}\right)  }\leq C\varepsilon.
\]

\section{The partial asymptotic domain decomposition}

The constructed above asymptotic expansions of the solutions of the Stokes
problem (\ref{problin}), (\ref{1.2010}), (\ref{1.2030}) and the
diffusion-convection-sorption problem (\ref{1.1000}), (\ref{1.2020}),
(\ref{1.2040}) allow us to apply the idea of the partial asymptotic domain
decomposition \cite{GP98, GP99}. We will cut off the two-dimensional
subdomains of $B^{\varepsilon}$ containing the "stenosis areas", bifurcations
and eventually (for the Stokes problem) the entrance and exit elements by the
lines orthogonal to the rectangles $B_{\varepsilon j}$ at the distance
$\delta=K\varepsilon\left\vert \ln\varepsilon\right\vert $ from the nodes
(bifurcation points) and from the nodal points of the stenosis areas. We will
call these subdomains the 2D zoom zones. Here $K$ is independent of
$\varepsilon$ and will be defined later.

Then we pass to the 1D description out of these subdomains. It means that we
pass to the projection of the variational formulation of the Stokes problem on
the Sobolev subspace of vector-valued functions having the Poiseuille
"parabolic" shape out of these subdomains (see \cite{BGPZ, GP05}).

For the diffusion-convection-sorption problem we apply the projection on the
Sobolev space of functions having vanishing derivatives of order greater than
$2k$ in the direction orthogonal to $B_{\varepsilon j}$ (also out of these 2D
zoom zones). This choice of the projection space is motivated by Lemma
\ref{lem1} and formula (\ref{cd3}). This gives us a variational formulation of
the partially decomposed diffusion-convection problem and according to the
general theory of the error estimate for the method of partial asymptotic
domain decomposition (\cite{GP05, GP99}) we get the estimates:

- for the difference of $u_{\varepsilon}$ and $u_{\varepsilon,\delta}^{dec}$
(solution of the partially decomposed problem for Stokes equation) we get as
in \cite{BGPZ} that for any $N$ there exists $K$ independent of $\varepsilon$
such that, if $\delta=K\varepsilon\left\vert \ln\varepsilon\right\vert $ then%
\[
\left\Vert u_{\varepsilon}-u_{\varepsilon,\delta}^{dec}\right\Vert
_{H^{1}\left(  B_{\varepsilon}\right)  }\leq O\left(  \varepsilon^{N}\right)
\sqrt{\varepsilon};
\]

- for the difference of $c_{\varepsilon}$ and $c_{\varepsilon,\delta}^{dec}$
(solution of the partially decomposed problem for the diffusion-convection
equation) we get that there exists $K$ independent of $\varepsilon$ such that,
if $\delta=K\varepsilon\left\vert \ln\varepsilon\right\vert $ then%
\begin{equation}
\left\Vert c_{\varepsilon}-c_{\varepsilon,\delta}^{dec}\right\Vert
_{H^{1}\left(  B_{\varepsilon}\right)  }\leq O\left(  \varepsilon^{k}\right)
\sqrt{\varepsilon}. \label{2.300}%
\end{equation}
These estimates justify the application of the MAPDD.

\begin{remark}
The interface conditions between the 2D parts of the domain and the 1D parts
follows from the variational formulation for the partially decomposed problem
by integrating by parts.
\end{remark}

\begin{remark}
Although for $k=0$ estimate (\ref{2.300}) is not too precise, we will hold
below a numerical experiment comparing the difference between the exact
solution and the solution of partially decomposed problem in this simpliest
case, when the projection space consists of functions with vanishing first
transversal derivative out of the 2D zoom zones.
\end{remark}

\section{Numerical experiments}

Here we will compare the solution of the leading term 1D equation
(\ref{1.500})-(\ref{1.6100}) to the numerical FEM solution of the coupled 2D
flow-diffusion problem (\ref{problin}), (\ref{1.1000}), (\ref{1.2020}),
(\ref{1.2030}), (\ref{1.2040}) in a thin rectangle $\left(  0,1\right)
\times\left(  0,\varepsilon\right)  $. We will trace the boundary layer
zones.
\[%
%TCIMACRO{\FRAME{itbpFU}{2.4405in}{0.8069in}{0in}{\Qcb{Fig. 4: Thin rectangle
%("straight channel geometry")}}{}{figure4n.eps}%
%width 2.4405in;  height 0.8069in;  depth 0in;  original-width 3.3053in;
%original-height 1.0758in;  cropleft "0";  croptop "1";  cropright "1";
%cropbottom "0";  filename 'Figure4N.eps';file-properties "XNPEU";}} }%
%BeginExpansion
{\parbox[b]{2.4405in}{\begin{center}
\includegraphics[
height=0.8069in,
width=2.4405in
]%
{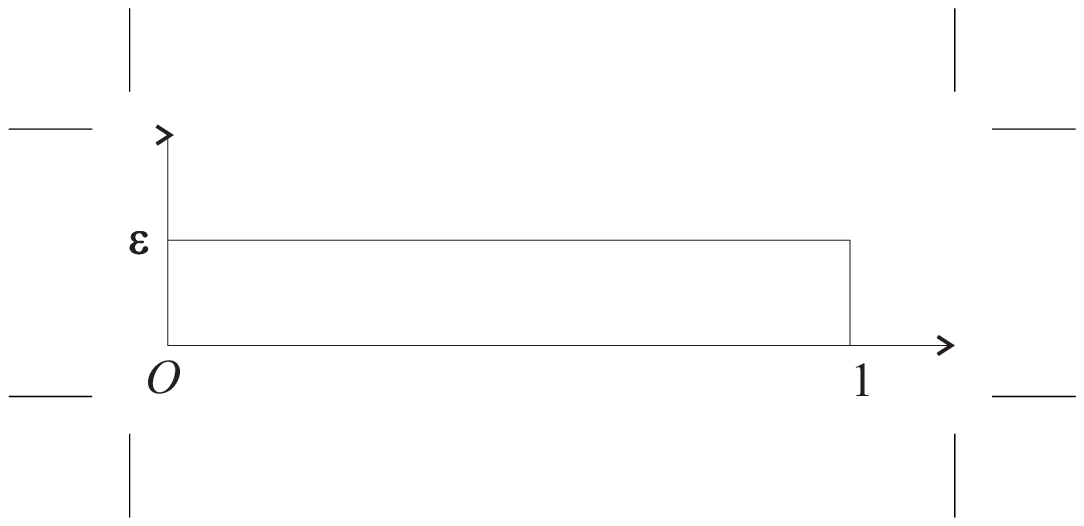}%
\\
Fig. 4: Thin rectangle ("straight channel geometry")
\end{center}}}
%EndExpansion
\]
In the second part of this section we will discuss the numerical solution
attained by the MAPDD for a one boundle tube structure corresponding to three
segments $e_{1},e_{2},e_{3}$.%
\[%
%TCIMACRO{\FRAME{itbpFU}{2.8677in}{2.3817in}{0in}{\Qcb{Fig. 5: One bundle tube
%structure}}{}{figure4.eps}{\special{ language "Scientific Word";
%type "GRAPHIC";  maintain-aspect-ratio TRUE;  display "USEDEF";
%valid_file "F";  width 2.8677in;  height 2.3817in;  depth 0in;
%original-width 4.9018in;  original-height 4.0629in;  cropleft "0";
%filename 'Figure4.eps';file-properties "XNPEU";}} }%
%BeginExpansion
{\parbox[b]{2.8677in}{\begin{center}
\includegraphics[
height=2.3817in,
width=2.8677in
]%
{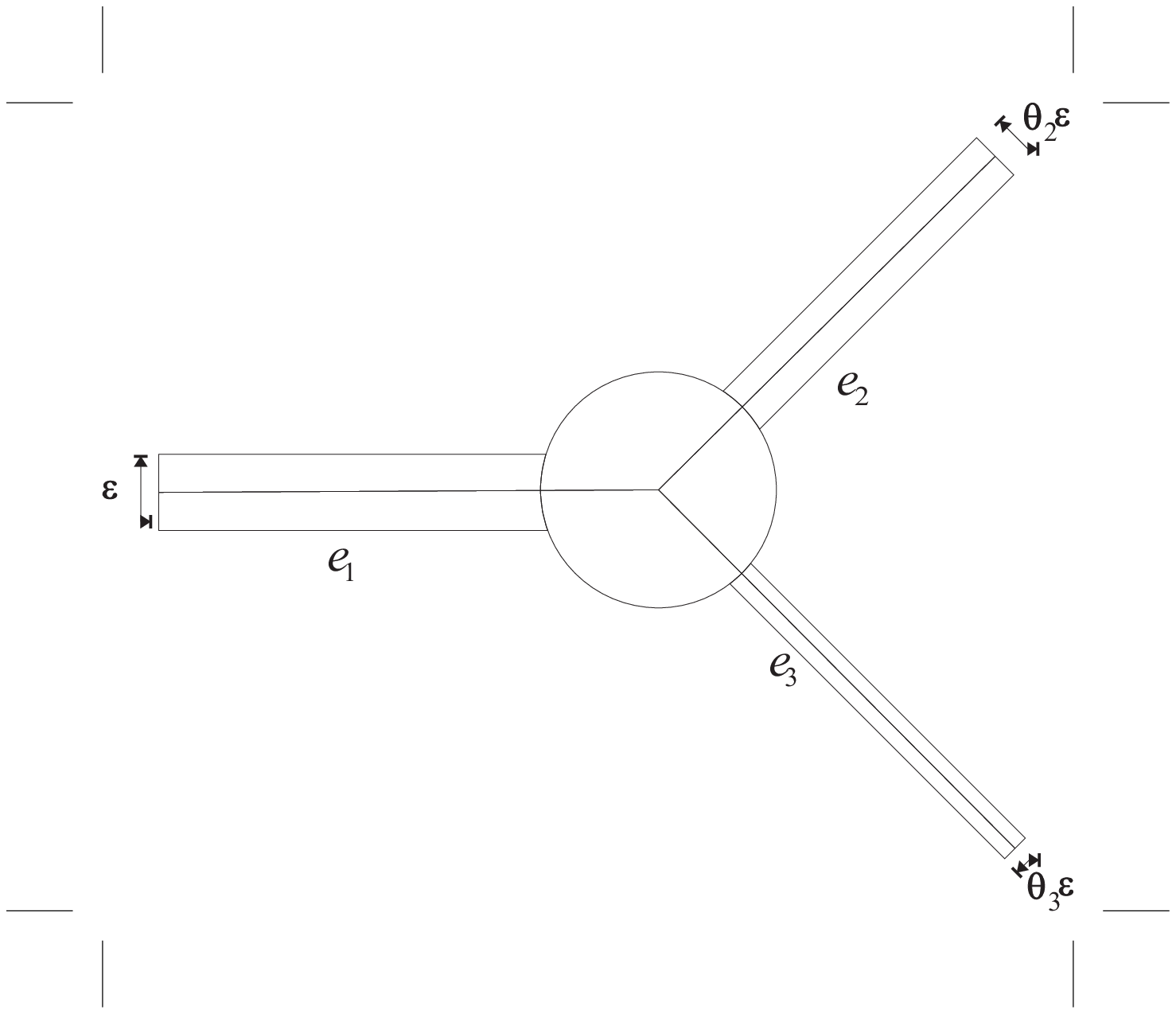}%
\\
Fig. 5: One bundle tube structure
\end{center}}}
%EndExpansion
\]

The finite element discretization of the Stokes equations (\ref{problin}) is
based on the classical P2/P1 lagrangian finite element test functions in
combination with the triangular finite element mesh. More precisely, the
velocity field is approximated by quadratic lagrangian test functions while
the pressure field is approximated with linear lagrangian test functions. As
it is well known, this finite element flow formulation satisfies the classical
Babuska-Brezzi \cite{C} compatibility condition and consequently produces
numerically stable and adequate solution strategy for Stokes and Navier-Stokes
problems. The concentration field $c$ is approximated by quadratic lagrangian
test functions. Finally, we would like to mention that for the studied flow
conditions, no specific divergence problems were encountered.

\subsection{Straight channel geometry}

Before analyzing the case of a 2D bifurcation problem, we will present some
results concerning the simple straight channel rectangle geometry
$B_{\varepsilon}=\left(  0,\varepsilon\right)  \times\left(  0,1\right)  .$
The main goal is to compare the predictions of our two methods (complete 2D
and asymptotic 1D model) in this simple flow conditions:$\ \varepsilon$ is
taken equal to $0.05,$ viscosity $\mu=1.$

In this case, the velocity distribution is described by a planar Poiseuille
flow given by formula (\ref{1.3000}). The entry concentration is fixed as
$q_{0}=1,$ while the exit concentration is maintained to be $q_{1}=0.5.$ We
calculate the solutions for $\beta=0.4$ and three different values of the
diffusion coefficient $\varkappa.$ The comparison of 2D solution and the 1D
asymptotic solution is presented at Fig. 6.

These variations of the diffusion coefficient are taken in order to find the
limits of the asymptotic approximation. This asymptotic analysis was applied
under the hypothesis of absence of other small parameters in the model.
Indeed, when the diffusion coefficient becomes a second small parameter then
the asymptotic analysis taking into account only one small parameter may be
not too precise.

$%
\begin{array}
[c]{c}%
%TCIMACRO{\FRAME{itbpF}{5.4483in}{1.4287in}{0in}{}{}{fig5an.eps}%
%width 5.4483in;  height 1.4287in;  depth 0in;  original-width 8.0575in;
%original-height 2.0825in;  cropleft "0";  croptop "1";  cropright "1";
%cropbottom "0";  filename 'Fig5aN.eps';file-properties "XNPEU";}} }%
%BeginExpansion
{\includegraphics[
height=1.4287in,
width=5.4483in
]%
{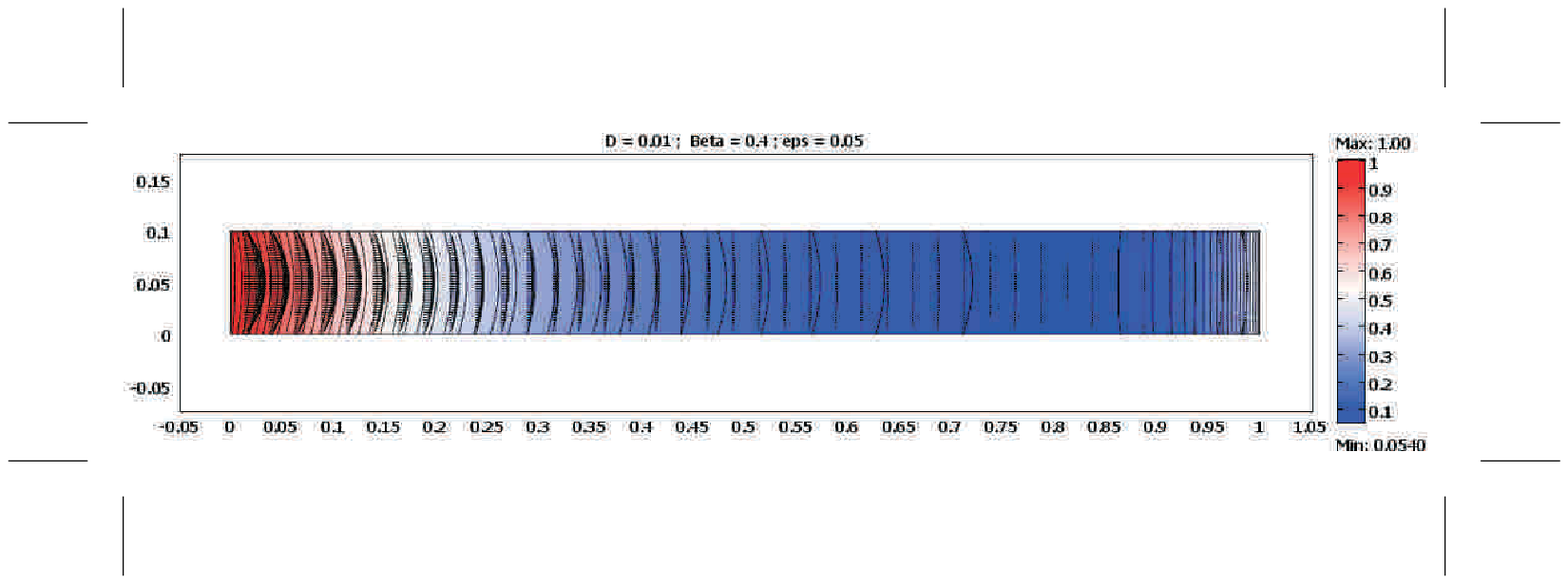}%
}
%EndExpansion
(a)\\%
%TCIMACRO{\FRAME{itbpF}{5.4474in}{1.4261in}{0in}{}{}{fig5bn.eps}%
%width 5.4474in;  height 1.4261in;  depth 0in;  original-width 8.0151in;
%original-height 2.0678in;  cropleft "0";  croptop "1";  cropright "1";
%cropbottom "0";  filename 'Fig5bN.eps';file-properties "XNPEU";}} }%
%BeginExpansion
{\includegraphics[
height=1.4261in,
width=5.4474in
]%
{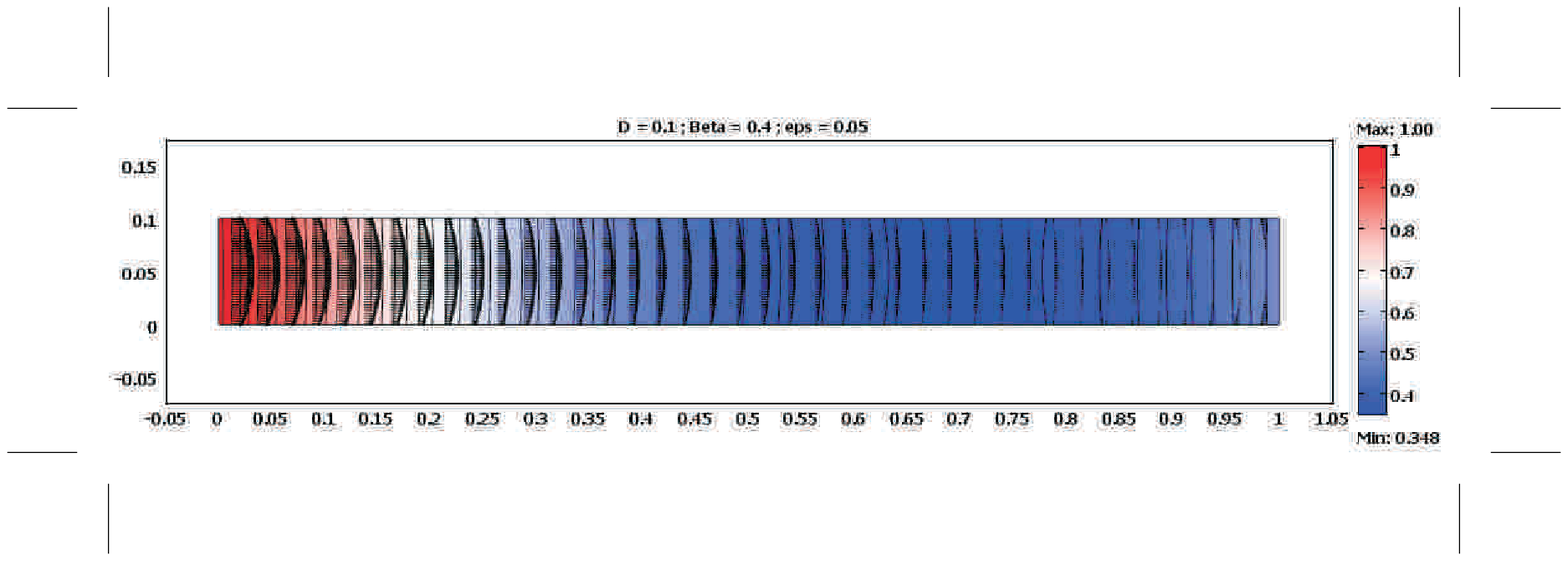}%
}
%EndExpansion
(b)\\%
%TCIMACRO{\FRAME{itbpF}{5.4483in}{1.4252in}{0in}{}{}{fig5cn.eps}%
%width 5.4483in;  height 1.4252in;  depth 0in;  original-width 8.1275in;
%original-height 2.0954in;  cropleft "0";  croptop "1";  cropright "1";
%cropbottom "0";  filename 'Fig5cN.eps';file-properties "XNPEU";}} }%
%BeginExpansion
{\includegraphics[
height=1.4252in,
width=5.4483in
]%
{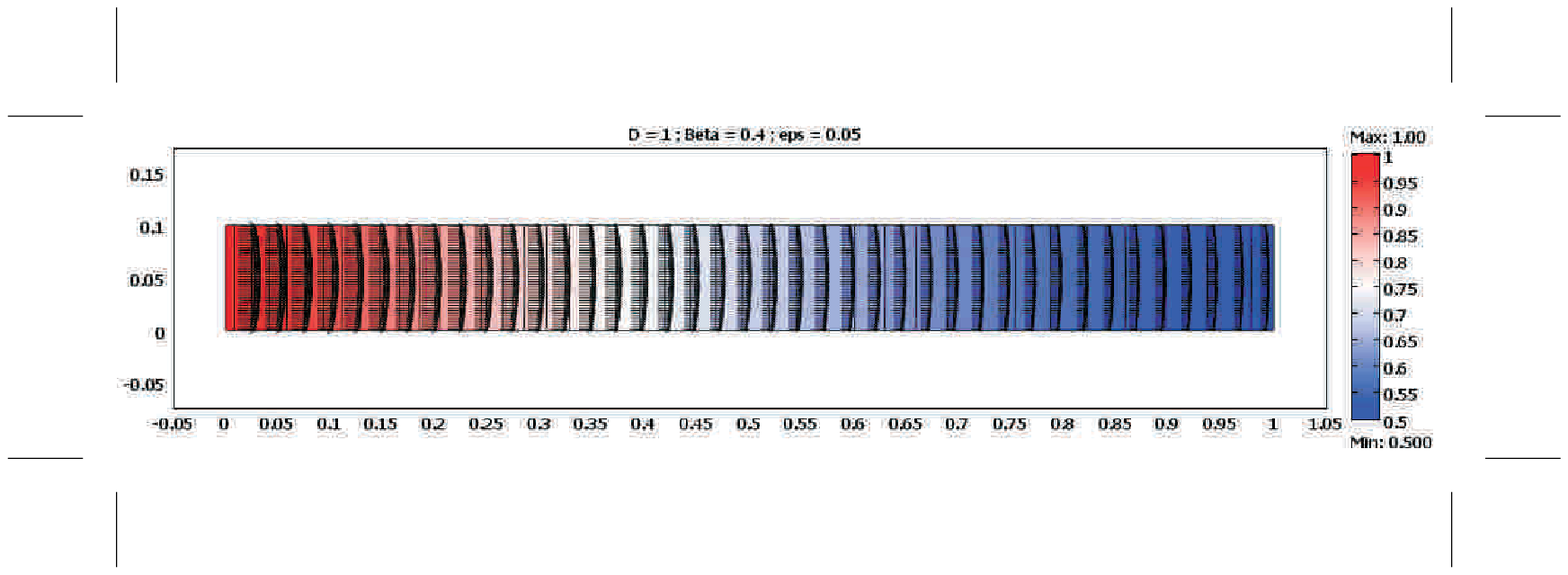}%
}
%EndExpansion
(c)\\
\text{Fig. 6: Distribution of the velocity field and distribution of the
concentration for values }\\
\text{(a) }\varkappa=0.01,\text{ (b) }\varkappa=0.1,\text{ (c) }\varkappa=1.\\
\text{Steady state concentration distribution: color map shows the
concentration distribution;}\\
\text{the arrows show that the total flux density for the concentration
distribution }\\
\text{(i.e., diffusive plus convective flux).}%
\end{array}
$

\bigskip

At low diffusion (Fig. 6(a) and 6(b)), the concentration distribution seems to
be essentially 2D and therefore one can expect some differences between the
corresponding 1D and 2D approaches. On the other hand, when the diffusion
increases, the concentration iso-levels become more and more planar and,
consequently, we can expect better performance from the 1D approach. This
analysis is further confirmed by the direct comparison between the predictions
of the 1D and the 2D approaches presented in Fig. 7.

$%
\begin{array}
[c]{c}%
%TCIMACRO{\FRAME{itbpF}{5.0842in}{3.4904in}{0in}{}{}{fig6n.eps}%
%width 5.0842in;  height 3.4904in;  depth 0in;  original-width 7.6268in;
%original-height 5.2157in;  cropleft "0";  croptop "1";  cropright "1";
%cropbottom "0";  filename 'Fig6n.eps';file-properties "XNPEU";}} }%
%BeginExpansion
{\includegraphics[
height=3.4904in,
width=5.0842in
]%
{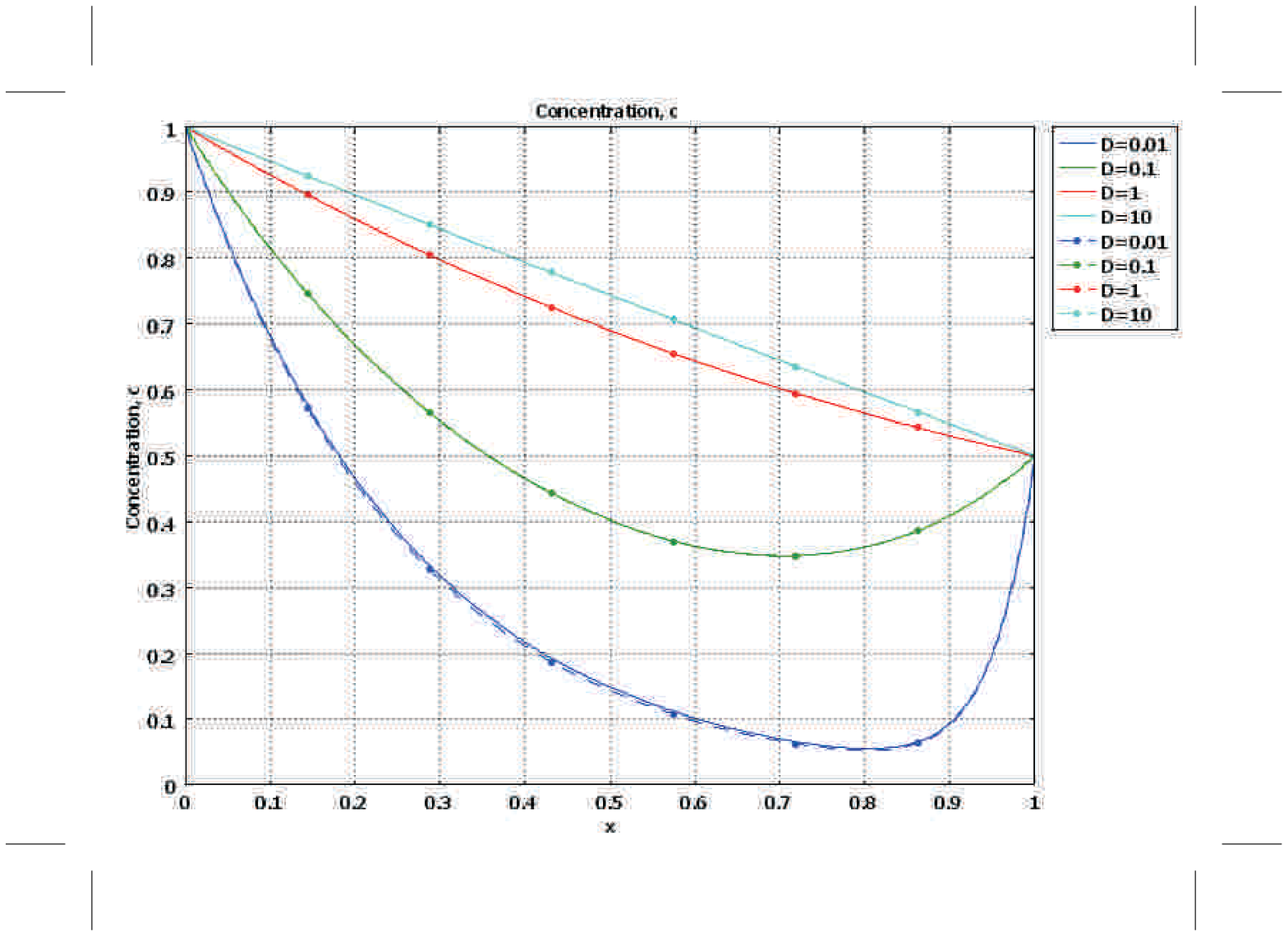}%
}
%EndExpansion
\\
\text{Fig. 7: Direct comparison between the 2D (full line) and 1D (dashed
line) numerical }\\
\text{predictions for 8 different values of the diffusion }\varkappa.\\
\text{Concentration distribution: full line stands for the average over a
cross-section of }\\
\text{concentration distribution for the 2D geometry; }\\
\text{dashed line and dots stand for the 1D concentration distribution.}%
\end{array}
$

\bigskip

As we can see from Fig. 7, some differences between 1D and 2D predictions
exist only for the value of $\varkappa=0.01$ when the diffusion is $5$ times
smaller than $\varepsilon.$ So the numerical experiment confirms the great
precision of the asymptotic solution (even in the case of small diffusion coefficient!).

\subsection{2D bifurcation geometry}

In this part, we extend our study to a more complex 2D bifurcation geometries.
The 2D flow geometry is presented in Fig. 8(a). As it is seen from this
figure, each channel could have different thickness and could be expressed in
terms of so called streamline function $\psi$ according to the following
definition:%
\begin{equation}
\frac{\partial\psi}{\partial x}=-U_{y};\ \frac{\partial\psi}{\partial y}=U_{x}
\label{1.800}%
\end{equation}
which is calculated as a solution of the following differential equation%
\begin{equation}
\Delta\psi=\frac{\partial U_{x}}{\partial y}-\frac{\partial U_{y}}{\partial
x}. \label{1.900}%
\end{equation}
The flow kinematics around the bifurcation point is illustrated in Fig. 8(a).
The corresponding pressure distribution is given in Fig. 8(b). As it is
predicted by the asymptotic analysis, the pressure gradient in each arm is
constant and naturally depends on the channel thickness and flow rate distribution.

$%
\begin{array}
[c]{cc}%
%TCIMACRO{\FRAME{itbpF}{2.9706in}{2.3826in}{0in}{}{}{fig8an.eps}%
%width 2.9706in;  height 2.3826in;  depth 0in;  original-width 6.2111in;
%original-height 4.9692in;  cropleft "0";  croptop "1";  cropright "1";
%cropbottom "0";  filename 'Fig8aN.eps';file-properties "XNPEU";}} }%
%BeginExpansion
{\includegraphics[
height=2.3826in,
width=2.9706in
]%
{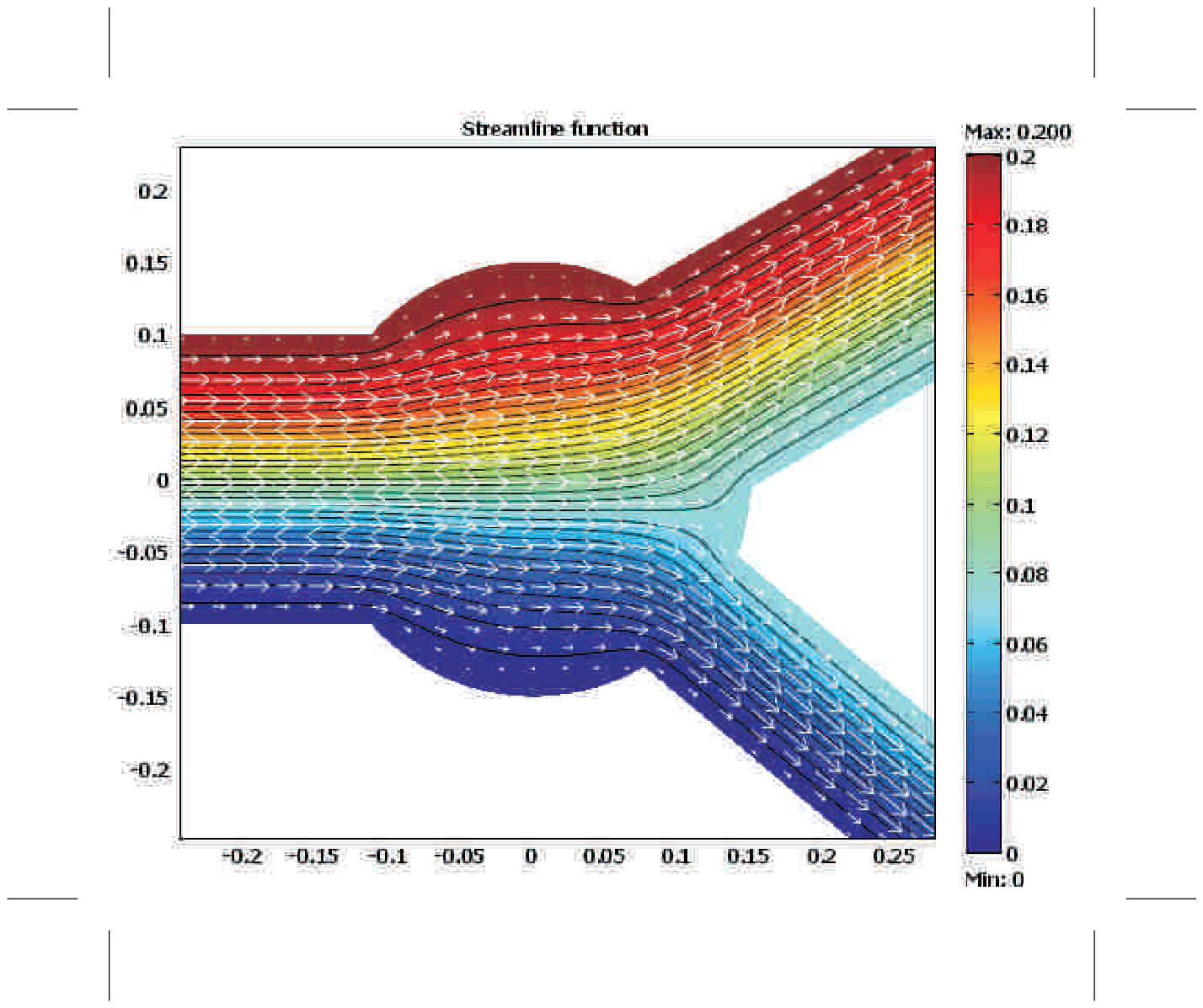}%
}
%EndExpansion
&
%TCIMACRO{\FRAME{itbpF}{2.9706in}{2.3826in}{0in}{}{}{fig8bn.eps}%
%width 2.9706in;  height 2.3826in;  depth 0in;  original-width 6.2111in;
%original-height 4.9692in;  cropleft "0";  croptop "1";  cropright "1";
%cropbottom "0";  filename 'Fig8bN.eps';file-properties "XNPEU";}} }%
%BeginExpansion
{\includegraphics[
height=2.3826in,
width=2.9706in
]%
{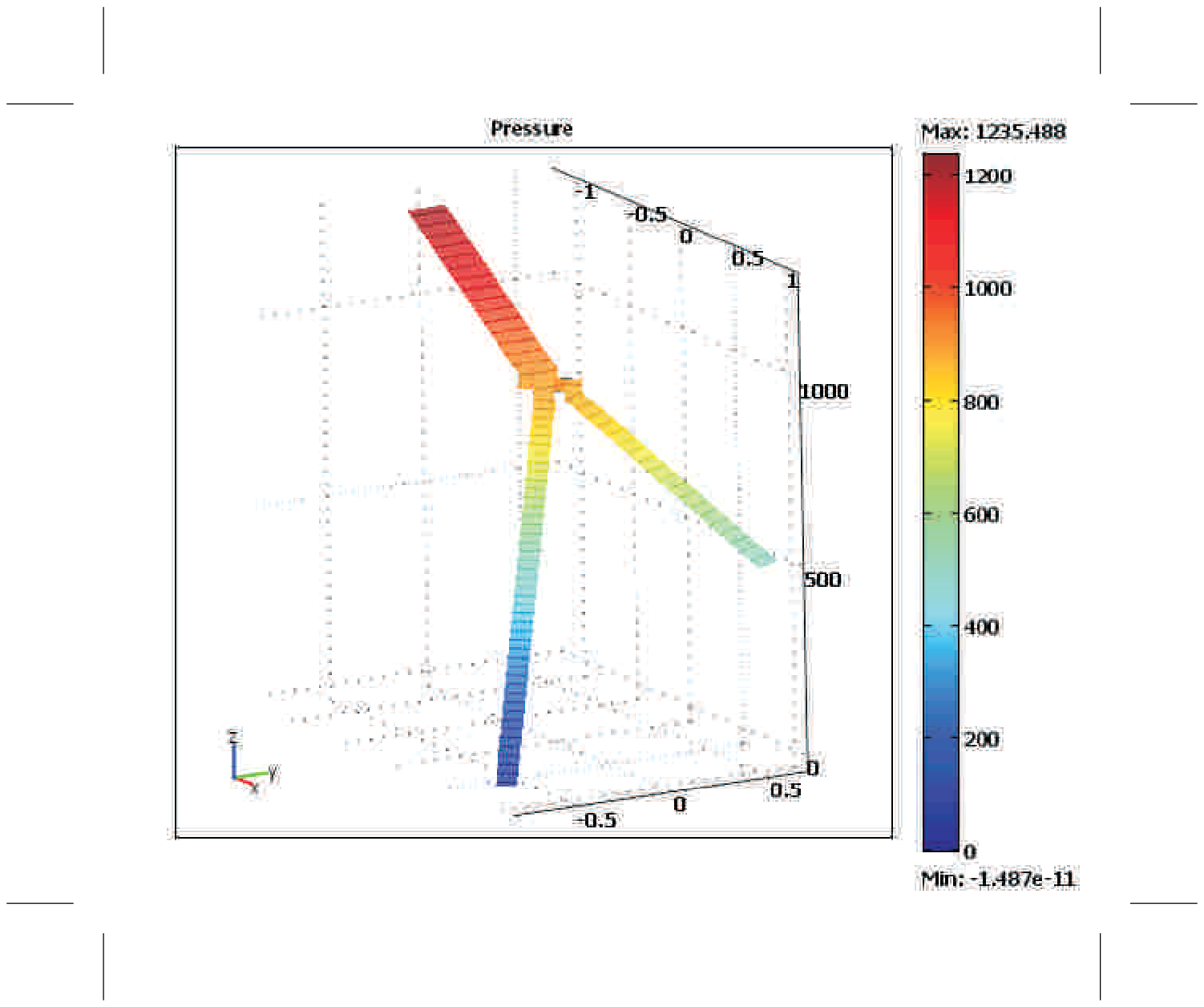}%
}
%EndExpansion
\\
(a) & (b)\\%
%TCIMACRO{\FRAME{itbpF}{3.3806in}{2.3817in}{0in}{}{}{fig8cn.eps}%
%width 3.3806in;  height 2.3817in;  depth 0in;  original-width 7.0707in;
%original-height 4.9649in;  cropleft "0";  croptop "1";  cropright "1";
%cropbottom "0";  filename 'Fig8cN.eps';file-properties "XNPEU";}} }%
%BeginExpansion
{\includegraphics[
height=2.3817in,
width=3.3806in
]%
{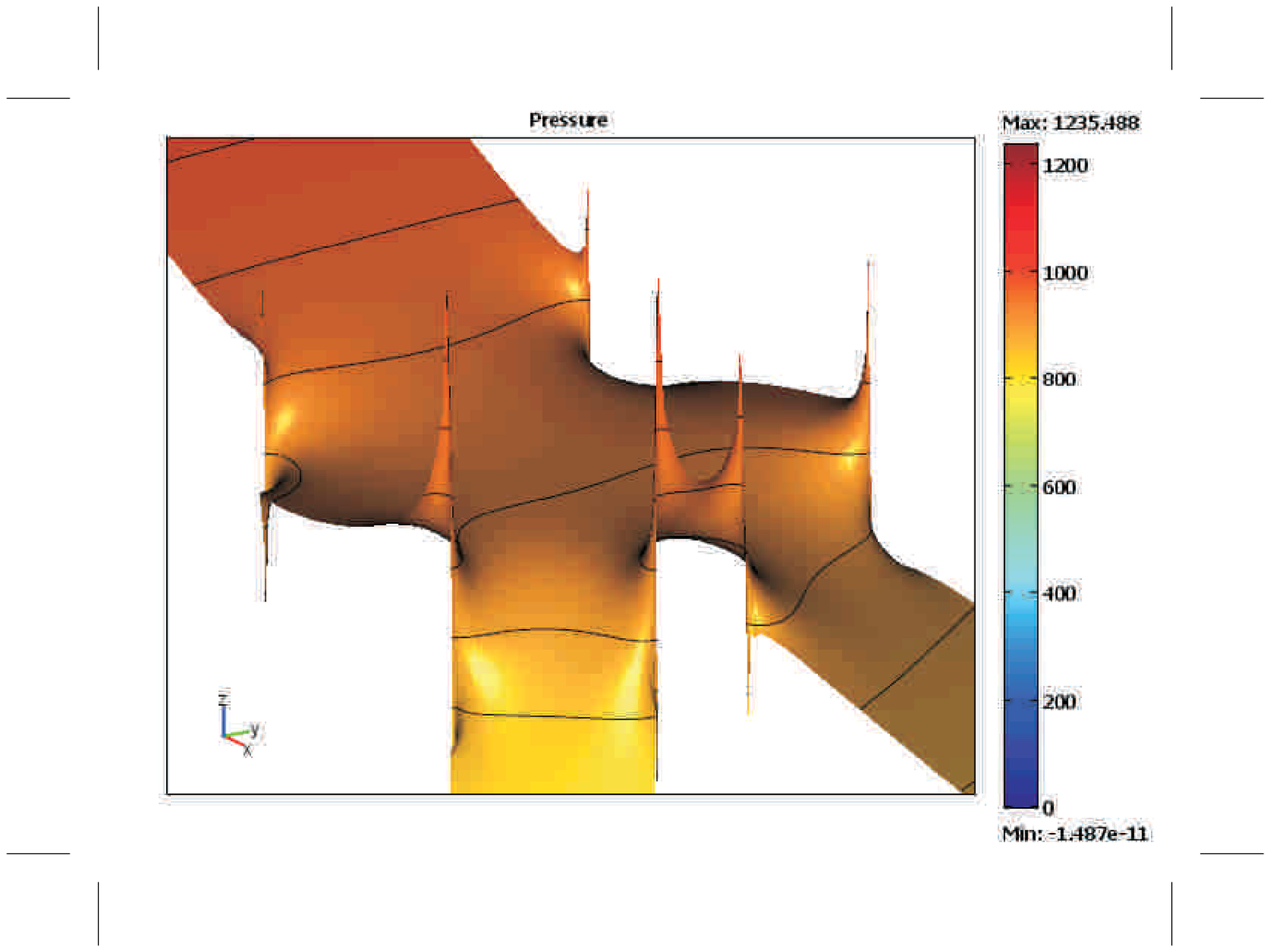}%
}
%EndExpansion
&
%TCIMACRO{\FRAME{itbpF}{3.3053in}{2.3817in}{0in}{}{}{fig8dn.eps}%
%width 3.3053in;  height 2.3817in;  depth 0in;  original-width 7.7652in;
%original-height 5.5763in;  cropleft "0";  croptop "1";  cropright "1";
%cropbottom "0";  filename 'Fig8dN.eps';file-properties "XNPEU";}} }%
%BeginExpansion
{\includegraphics[
height=2.3817in,
width=3.3053in
]%
{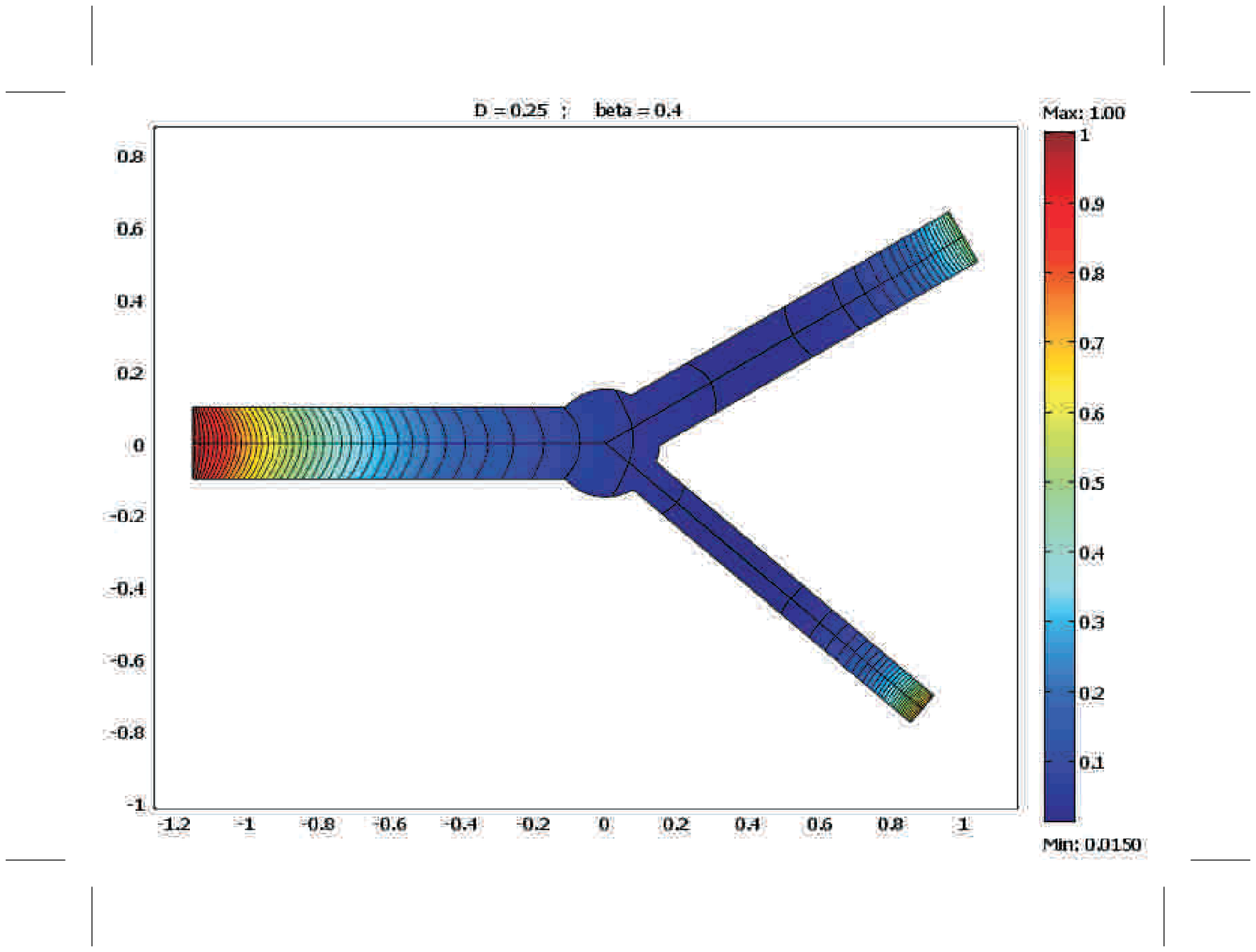}%
}
%EndExpansion
\\
(c) & (d)\\%
%TCIMACRO{\FRAME{itbpF}{3.3053in}{2.3817in}{0in}{}{}{fig8en.eps}%
%width 3.3053in;  height 2.3817in;  depth 0in;  original-width 7.7652in;
%original-height 5.5763in;  cropleft "0";  croptop "1";  cropright "1";
%cropbottom "0";  filename 'Fig8eN.eps';file-properties "XNPEU";}} }%
%BeginExpansion
{\includegraphics[
height=2.3817in,
width=3.3053in
]%
{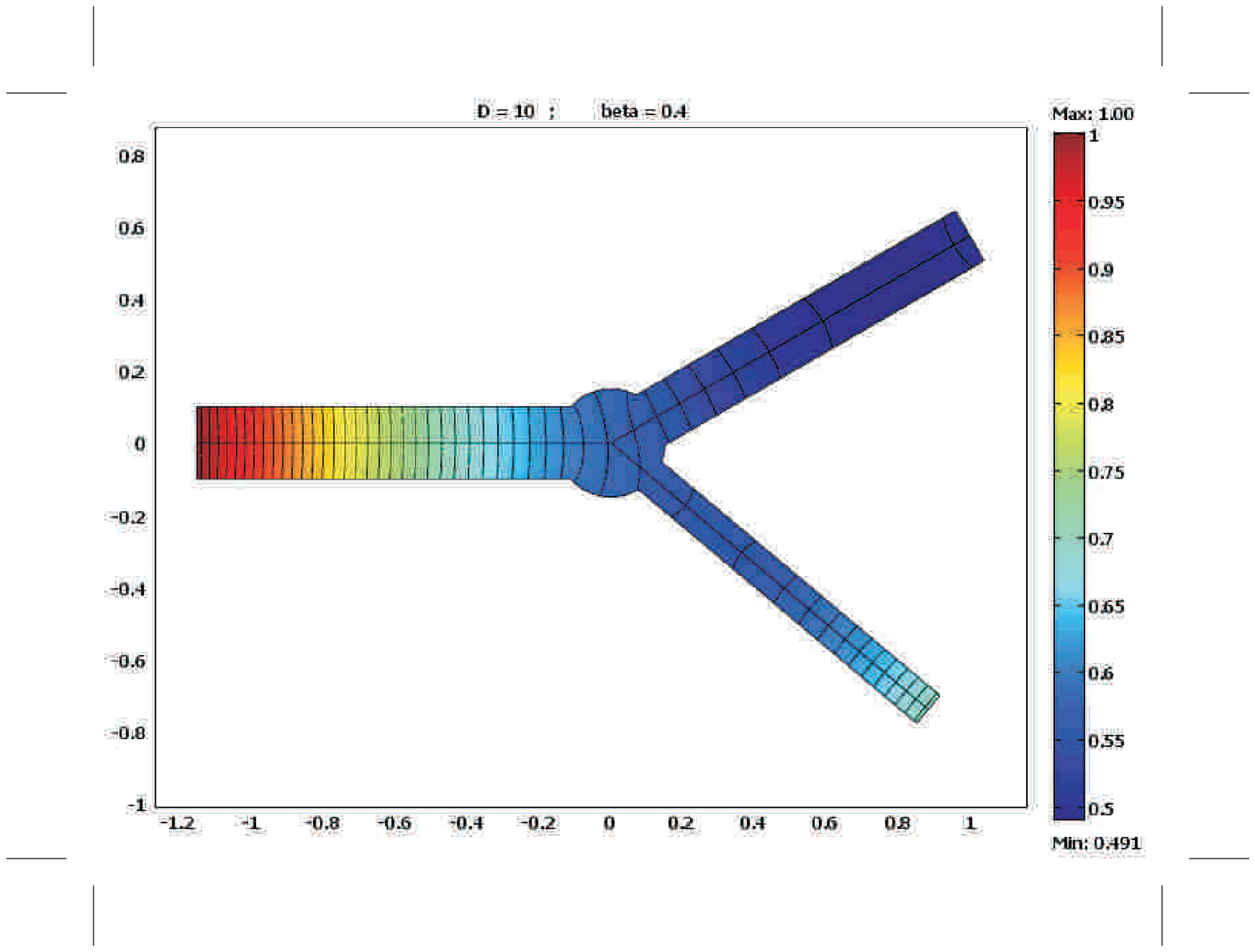}%
}
%EndExpansion
&
%TCIMACRO{\FRAME{itbpF}{3.704in}{2.3817in}{0in}{}{}{fig8fn.eps}%
%width 3.704in;  height 2.3817in;  depth 0in;  original-width 7.7513in;
%original-height 4.9649in;  cropleft "0";  croptop "1";  cropright "1";
%cropbottom "0";  filename 'Fig8fN.eps';file-properties "XNPEU";}} }%
%BeginExpansion
{\includegraphics[
height=2.3817in,
width=3.704in
]%
{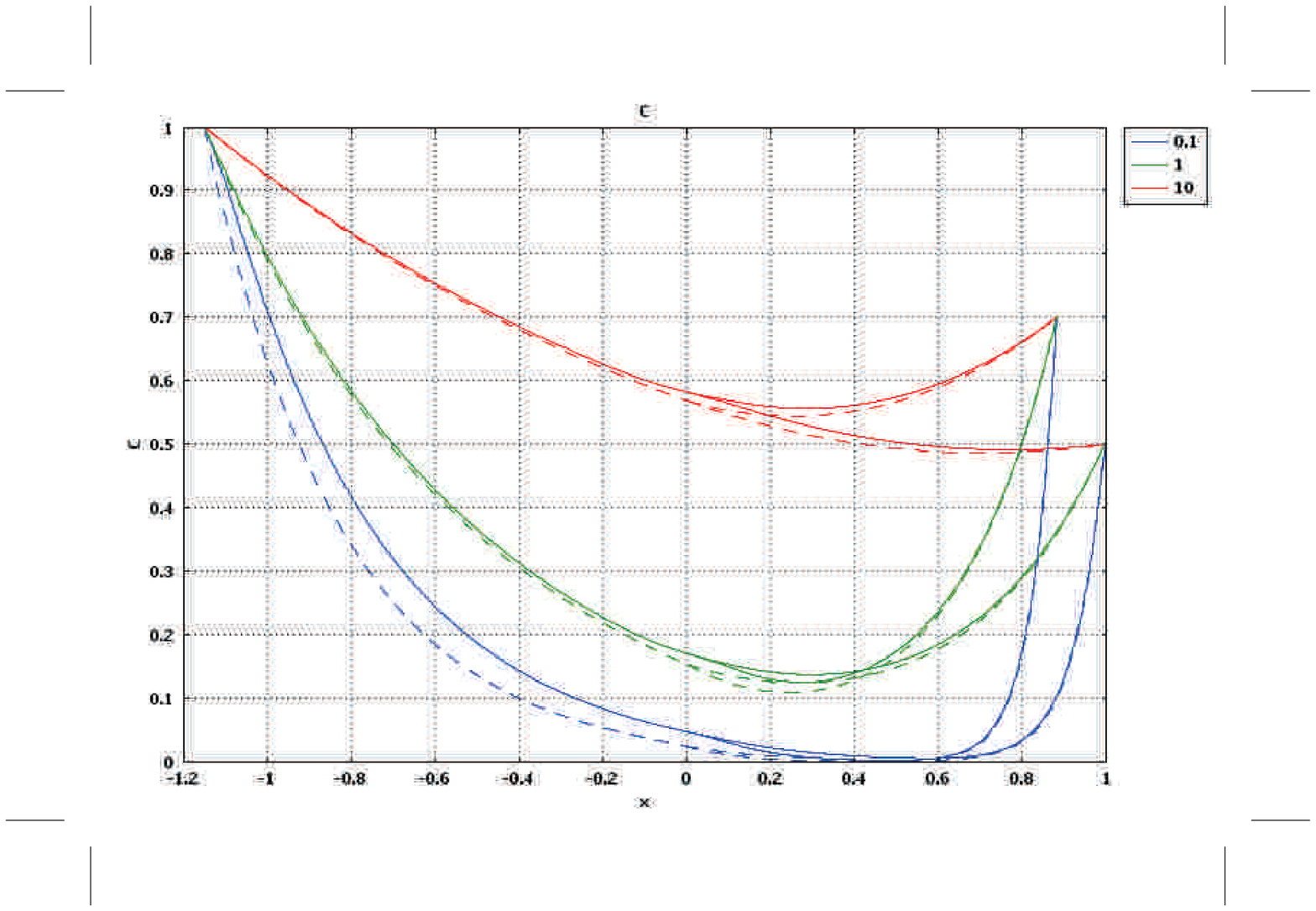}%
}
%EndExpansion
\\
(e) & (f)
\end{array}
$

Fig. 8: Model of the Stokes flow and the convection/diffusion process in a one
boundle tube structure: (a) flow kinematics, (b), (c) pressure distribution in
the bifurcation zone, (d) and (e) iso-lines of concentration for two values of
$\varkappa:$ $\varkappa=0.25$ (d) and $\varkappa=0.10$ (e), (f) comparison
between the 2D and MAPDD solutions for the diffusion $\varkappa=0.1,$ $1$ and
$10$.

\bigskip

It is important to emphasize that each geometry corner represents a singular
point for the pressure, which can be seen at Fig. 8(c). These peaks are
predicted by the corner singularities analysis (\cite{Naz, Grisvard, Nicaise,
Sol}). Two typical solutions representing the cases of lower ($\varkappa
=0.25$) and higher ($\varkappa=10$) diffusivity are given in Fig. 8(d) and
8(e). Like in the previously analyzed geometries, there is a critical
diffusivity value, which ensures the validity of the simplified 1D approach.
In Fig. 8(f), we have given the direct comparison between 1D and 2D
predictions. As it is seen from these figures, the critical diffusivity value
is around 1.


\begin{thebibliography}{99}                                                                                               %


\bibitem {GPbook}G.P. Panasenko, \textit{Multi-Scale Modelling for Structures
and Composites}, Springer, 2005.

\bibitem {GP91}Panasenko G.P., \textit{Asymptotic solutions of the elasticity
theory system of equations for lattice and skeletal structures, }Math. Sb.,
183 (1), 1992, 89-113; Engl. Transl. in Russian Acad. Sci. Sbornik Math. 75
(1), 1993, 85-110.

\bibitem {GP00}Panasenko G.P., \textit{Asymptotic expansion of the solution of
Navier-Stokes equation in tube structure and partial asymptotic decomposition
of the domain}, Appl. Anal. 76 (2000), no. 3-4, 363--381.

\bibitem {Naz}Nazarov S.A., Plamenevskii, \textit{Elliptic problems in domains
with piecewise smooth boundaries}, Berlin-New York Walter de Gruyter, 1994.

\bibitem {Galdi}Galdi G.P., \textit{An introduction to the mathematical theory
of the Navier-Stokes equations I}, Springer-Verlag, 1994.

\bibitem {Lad}Ladyzhenskaya O.A., \textit{Boundary value problems of
Mathematical Physics}, Springer-Verlag, 1985.

\bibitem {GP98}Panasenko G.P., \textit{Method of asymptotic partial
decomposition of domain}, Mathematical Models and Methods in Applied Sciences,
8, n.1, 1998, 139-156.

\bibitem {BGPZ}Blanc F., Gipouloux O., Panasenko G.P., Zine A. M.,
\textit{Asymptotic analysis and partial asymptotic decomposition of domain for
Stokes equation in tube structure}, Math. Models Methods Appl. Sci. 9 (1999),
no. 9, 1351--1378.

\bibitem {GP05}Meliani S., Panasenko G.P., \textit{Thermo-chemical modelling
of the formation of a composite material}, Appl. Anal. 84 (2005), no. 3, 229--245.

\bibitem {GP99}Panasenko G.P., \textit{Asymptotic partial decomposition of
variational problems, C.R. Mechanique Acad. Sci. Paris, Serie IIb, 327, 1999,
1185-1190.}

\bibitem {C}Ciarlet P.G., \textit{The finite element method for elliptic
problems}, North-Holland, Amsterdam, 1978.

\bibitem {PSS}Panasenko G.P., Sirakov I., Stavre R., \textit{Asymptotic and
numerical modeling of a flow in a thin channel with viscoelastic wall}, Int.
Journal for Multiscale Computational Engineering, 5 (6), 2007.

\bibitem {Grisvard}Grisvard P., \textit{Elliptic Problems in Nonsmooth
Domains}, Pitman, 1985.

\bibitem {Nicaise}Nicaise S., \textit{About the Lam\'{e} system in a polygonal
or a polyhedral domain and a coupled problem between the Lame system and the
plate equation I: Regularity of the solutions,} Annali della Scuola Normale
Sup. di Pisa, ser. IV. 20, 1993, 163-191.

\bibitem {Sol}Solonnikov V. A., \textit{On the Stokes equations in domains
with non-smooth boundaries and on viscous incompressible flow with a free
surface}, College de France Seminar, 4, 1983, 240-349.
\end{thebibliography}
\end{document}